\documentclass{amsart}

\usepackage{color}
\usepackage{amssymb, amsthm, amsmath, amscd}
\usepackage{hyperref}

\usepackage[pdftex]{graphicx}
\usepackage{subfig}
\usepackage{float}
\usepackage[below]{placeins}

\newtheorem{lemma}{Lemma}
\newtheorem{proposition}[lemma]{Proposition}
\newtheorem{theorem}[lemma]{Theorem}
\newtheorem*{vartheorem}{Theorem 1}

\newtheorem{definition}[lemma]{Definition}
\newtheorem{corollary}[lemma]{Corollary}

\newcommand{\pd}{ {\partial}}
\newcommand{\R}{\mathbf R}
\newcommand{\eps}{\varepsilon}
\newcommand{\ubar} {\underline}

\newcommand{\lin}{\mathcal{L}}

\newcommand{\F}{\mathcal{F}}
\newcommand{\hcal}{\mathcal{H}}
\newcommand{\ncal}{\mathcal{N}}
\newcommand{\pcal}{\mathcal{P}}

\newcommand{\ical}{\mathcal{I}}

\newcommand{\dcal}{\mathcal{D}}

\newcommand{\ccal}{\mathcal{C}}

\begin{document}

\title[Self-Similar Surfaces under MCF III]{Construction of Complete Embedded Self-Similar Surfaces under Mean Curvature Flow. Part III. }
\author{Xuan Hien Nguyen}
\thanks{This work is partially supported by the National Science Foundation, Grant No. DMS-0908835.}
\address{Department of Mathematics, Iowa State University, Ames, IA 50011} 
\email{xhnguyen@iastate.edu}


\begin{abstract}
We present new examples of complete embedded self-similar surfaces under mean curvature flow by gluing a sphere and a plane. These surfaces have finite genus and are the first examples of non-rotationally symmetric self-shrinkers in $\R^3$. Although our initial approximating surfaces are asymptotic to a plane at infinity, the constructed self-similar surfaces are asymptotic to cones at infinity. 
\end{abstract}

\maketitle

\section{introduction}

This article is the third and last installment of a series of papers aiming at constructing new examples of surfaces satisfying the self-shrinking equation for the mean curvature flow, 
	\begin{equation}
	\label{eq:self-shrinker}
	\tilde H + \tilde X \cdot \nu =0,
	\end{equation}
where $\tilde X$ is the position vector, the function $\tilde H$ and the orientation of the unit normal $\nu$ are taken so that the mean curvature vector is given by $\tilde {\mathbf H} = \tilde  H \nu$. 

In \cite{huisken;asymptotic-behavior}, Huisken proved that if the growth of the second fundamental form $|A|^2$ is controlled (type 1), the singularities of the mean curvature flow tend asymptotically to a solution to \eqref{eq:self-shrinker}. The work on self-shrinking surfaces is therefore motivated by a desire to better understand the regularity of the mean curvature flow. A long list of examples of self-shrinkers would help shed light on the behavior of the flow near its singularities; unfortunately, until now, there were only four known examples of complete embedded self-shrinking surfaces (in the Euclidean space $E^3$): a plane, a cylinder, a sphere, and a shrinking doughnut \cite{angenent;doughnuts}; although there is numerical evidence of many others   \cite{angenent-chopp-ilmanen;computed-example-nonuniqueness}\cite{chopp;computation-solutions}.

The overarching idea in the three articles is to obtain new examples of self-shrinkers by desingularizing the intersection of two known examples  (the sphere of radius $\sqrt 2$ centered at the origin and a plane through the origin) using Scherk minimal surfaces. First, one constructs an initial approximate solution $\tilde M$ by fitting an appropriately bent and scaled Scherk surface $\tilde \Sigma$ in a neighborhood of the intersection, then one solves a perturbation problem in order to find an exact solution. The method was successfully used by Kapouleas \cite{kapouleas;embedded-minimal-surfaces} and Traizet \cite{traizet;surfaces-minimales} to construct minimal surfaces, and by the author for self-translating surfaces under the mean curvature flow \cite{mine;tridents} \cite{mine;self-trans}.

The main difficulty in these desingularizations lies in showing that the linearized equation $\tilde \lin v := \Delta v + |\tilde A|^2 v - \tilde X \cdot \nabla v + v=E$  can be solved on the initial surface $\tilde M$. One attacks it by studying $\tilde \lin v =E$ on smaller pieces first. In the first article \cite{mine;part1}, we study the linearized equation on the desingularizing surface $\tilde \Sigma$. The second article concerns the outer plane $\tilde \pcal$ (the plane with a central disk removed) and its main result states that the Dirichlet problem for \eqref{eq:self-shrinker} on $\tilde \pcal$ has a unique solution among graphs of functions over $\tilde \pcal$ with a controlled linear growth. In the present article, we finish the construction by gluing the solutions to the linearized equations on the different pieces to obtain a global solution with a standard method: we use cut-off functions to localize the inhomogeneous term to the different pieces, solve the linearized equation on these pieces, glue the local solutions using cut-off functions again, and iterate the process. However, the  cut-off functions create errors and obtaining the right estimates for  the iteration to converge requires a delicate and precise construction of the initial approximate surfaces, which is the main focus of this article. 

Once the initial surfaces are constructed, the techniques 
from \cite{kapouleas;embedded-minimal-surfaces} can be readily applied, with one notable exception. In all the previous constructions \cite{kapouleas;embedded-minimal-surfaces}, \cite{traizet;surfaces-minimales}, \cite{mine;tridents}, and \cite{mine;self-trans}, the surfaces converge exponentially to their asymptotic catenoids, planes, or grim reapers respectively. But here, the self-shrinkers grow linearly at infinity. In this article, we also  refine previous estimates from \cite{mine;part2} and define the appropriate Banach spaces of functions for the final fixed point theorem. 

The author originally intended to desingularize the intersection of a cylinder and a plane, inspired by the  work of Kapouleas on minimal surfaces and Ilmanen's open problem on the rigidity of the cylinder. In \cite{ilmanen;lectures-mcf}, Ilmanen conjectured that if one of the ends of a self-shrinking surface is asymptotic to a cylinder, then the surface must be the self-shrinking cylinder itself. At the time, Angenent suggested that the desingularization of a cylinder and a plane could provide a counterexample. The problem was harder than anticipated because of the asymptotic behavior along the cylindrical ends and remains open. In the same set of lecture notes, Ilmanen predicted the existence of self-shrinking punctured saddles, which he called $m$-saddles, and provided pictures for $m=1, 3,$ and $9$. Theorem \ref{thm:8-2} below shows that these $m$-saddles exist if $m$ is large; the existence for  $m$ small is still an open problem. The parameter $m$ is the number of periods of a Scherk surface needed to wrap around the intersection circle. The $m$-saddle $\tilde M_m$ below therefore has $2m$ handles and genus $m-1$.

\begin{theorem} 
\label{thm:8-2}

There exists a natural number $\bar  m$ so that for any natural number $m>\bar m$, there exists a surface $\tilde M_m$ with the following properties:
	\begin{enumerate}
	\item $\tilde M_m$ is a complete smooth surface which satisfies the equation $\tilde H + \tilde X \cdot \nu =0$.
	\item $\tilde M_m$ is invariant under rotation of $180^{\circ}$ around the $\tilde x$-axis.
	\item $\tilde M_m$ is invariant under reflections across planes containing the $\tilde z$-axis and forming angles $\pi/(2m) + k \pi/m$, $k \in \mathbf Z$, with the $\tilde x$-axis.
	\item Let $U =B_2 \cap \{\tilde z>0\}$ be the open top half of the ball of radius $2$. As $m \to \infty$, the sequence of surfaces $\tilde M_m$ tends to the sphere of radius $\sqrt 2$ centered at the origin on any compact set of  $U$. 
	\item $\tilde M_m$ is asymptotic to a cone.
	\item If we denote by $T$ the translation by the vector $-\sqrt 2 \vec e_x$, the sequence of surfaces $m T(\tilde M_m) = \{ (m\tilde x, m\tilde y, m\tilde z) \mid (\tilde x + \sqrt 2 \vec e_x, \tilde y, \tilde z) \in \tilde M_m \} $  converges in $C^k$ to the original Scherk surface $\Sigma_0$ on compact sets, for all $k \in \mathbf N$.
	\end{enumerate}

\end{theorem}

We briefly sketch the proof below, highlighting the differences and similarities between this construction and the ones from \cite{kapouleas;embedded-minimal-surfaces} and  \cite{mine;self-trans}. 

We start by replacing a small neighborhood of the intersection circle by an appropriately bent Scherk surface to obtain embedded surfaces.  However, instead of scaling down the Scherk surface by a factor $ \tau$ where $\tau$ is a small positive constant, we keep it in its natural scale so that the curvatures and second fundamental form stay bounded, and scale up the rest of the configuration by $ \tau^{-1}$. The equation to be satisfied is then
	\begin{equation}
	\label{eq:self-shrinker-large}
	H + \tau^2 X \cdot \nu=0.
	\end{equation}
These initial surfaces are embedded and will be our approximate solutions. The next and more difficult step consists in finding an exact solution among perturbations of the initial surfaces. More precisely, we perturb a surface by adding the graph of a small function $v$ so the position vector $X$ becomes $X_v:=X+v \nu$.  Denoting the initial surface by $M$, its position vector by $X$, its unit normal vector by $\nu$, the graph of $v$ over $M$ by $M_v$, its mean curvature by $H_v$, and its unit normal vector by $\nu_v$, we have
	\[
	H_v + \tau^2 X_v \cdot \nu_v = H +\tau^2 X \cdot \nu + \Delta v + |A|^2 v - \tau^2 X \cdot \nabla v  + \tau^2 v + Q_v,
	\]
where  $A$ is the second fundamental form on $M$ and $Q_v$ is at least quadratic in $v$, $\nabla v$, and $\nabla^2 v$.  The surface $M_v$ is a self-shrinker if 
	\begin{equation}
	\label{eq:interm}
	\lin v = -H -\tau^2  X \cdot \nu - Q_v,
	\end{equation}
where $\lin v :=  \Delta v + |A|^2 v - \tau^2  X\cdot \nabla v+ \tau^2 v$. Once we can solve the equation  $\lin v =  -H -\tau^2  X \cdot \nu$, we expect the quadratic term to be small  so the solution $v$ to \eqref{eq:interm} could be obtained by iteration. Before we can solve the linearized equation $\lin v =E$ on the initial surface $M$, we have to study its associated Dirichlet problem on the various pieces: the desingularizing surface $\Sigma$ (formed by a truncated bent Scherk surface), the two rotationally symmetric caps $\mathcal{C}$, the inner disk $\mathcal{D}$, and the outer plane $\mathcal{P}$.  

In all of the previous constructions (and here also), the linear operator $\lin$  has small eigenvalues on $\Sigma$. One way to deal with the presence of small eigenvalues is to restrict the class of possible perturbations and eigenfunctions by imposing symmetries on all the surfaces considered. However, this method only works if the initial configuration has the imposed symmetries and, in general, can not rule out all  the troublesome eigenfunctions. A second complementary approach is to invert the linear operator modulo the eigenfunctions corresponding to small or vanishing eigenvalues. In other words, one can add or subtract a linear combination of  eigenfunctions to the inhomogeneous term of $\lin v=E$ in order to land in the space perpendicular to the approximate kernel, where the operator has a bounded inverse. For an exact solution, one must be able to generate (or cancel) any linear combination of these eigenfunctions within the construction. The process is called unbalancing and consists in dislocating the Scherk surface so that opposite asymptotic planes are no longer parallel. Flexibility in the initial configuration is the key to a successful construction. The terms unbalancing and flexibility were first introduced by Kapouleas and the reader can find a short discussion of them in the survey article \cite{kapouleas;clay}.  

\subsection{How this construction differs from previous ones.}

In \cite{kapouleas;embedded-minimal-surfaces} (\cite{mine;self-trans}), the flexibility relies on the fact that the main equation $\tilde H=0$  ($\tilde H- \vec e_y \cdot \nu=0$ resp.) is translation invariant, so the catenoidal ends  (grim reaper ends  resp.) could be shifted without creating errors. Moreover, since catenoids (grim reapers resp.) have ends,  one can, with careful planning, perform the required dislocation at the intersections so that  all the small changes in position build up toward loose ends. For the case of self-shrinkers, the sphere of radius $\sqrt 2$ centered at the origin is the only sphere satisfying \eqref{eq:self-shrinker} so the apparent lack of flexibility has been a major impediment in completing the desingularization of the sphere and the plane.

\begin{figure}[h]
\includegraphics[height=1in]{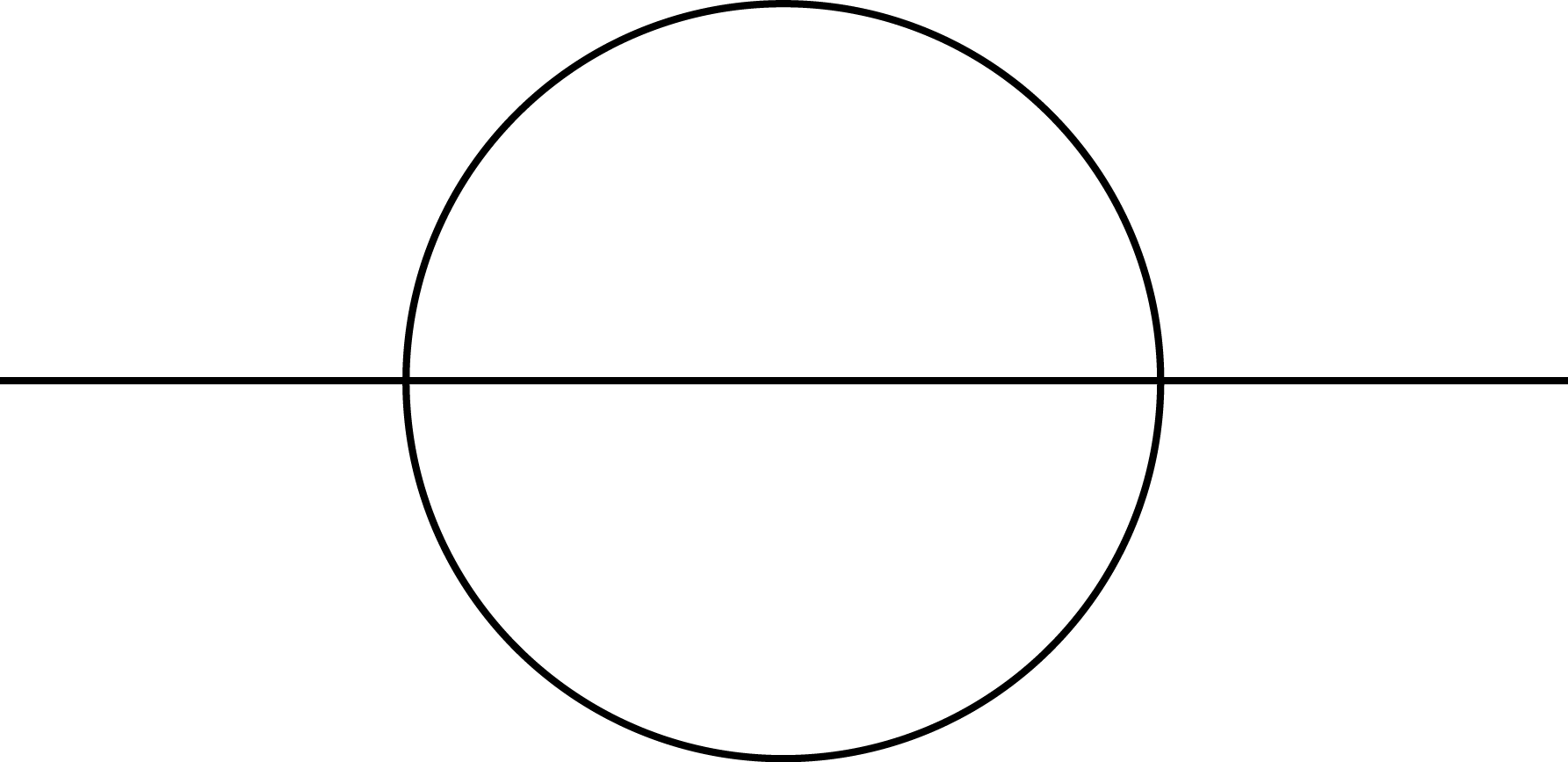}
\includegraphics[height=1in]{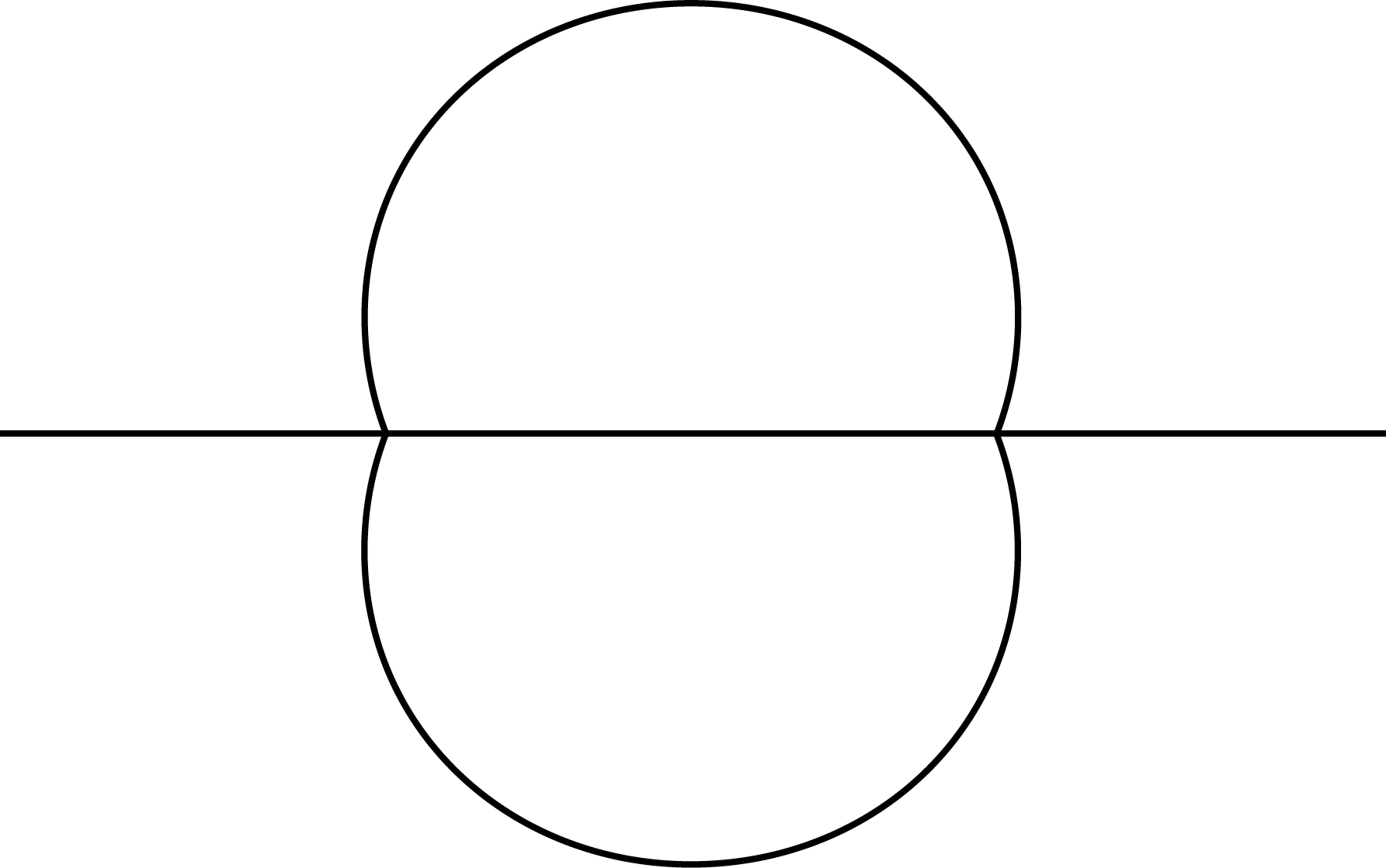}
\caption{A balanced initial configuration and an unbalanced one.}
\label{fig:simple-conf}
\end{figure}

The \emph{unbalancing} process requires one to consider the configuration of a sphere and a plane as part of a family of initial configurations in which the rotationally symmetric caps meet the plane at various angles close to $90$ degrees (see Figure \ref{fig:simple-conf} for a dramatized representation).

%

Rather than shifting the sphere up or down, which would create too much error, we use a family of  self-shrinking rotationally symmetric caps. In \cite{angenent;doughnuts}, Angenent showed that rotationally symmetric self-shrinkers are generated by geodesics in the half-plane $\{ (\tilde z, \tilde r) \mid \tilde r \geq 0\}$ with metric 
	\begin{equation}
	\label{eq:metric}
	\tilde r^{2} e^{-(\tilde z^2+\tilde r^2)} \{ (d\tilde z)^2 +(d\tilde r)^2\}.
	\end{equation}
The equation for these geodesics parametrized by arc length is given by the following system of ordinary differential equations:
	\begin{equation}
	\label{eq:system-ode}
	\begin{cases}
	\dot {\tilde z} = \cos \theta, \\
	\dot {\tilde r} = \sin \theta, \\
	\dot \theta = \tilde z \sin \theta + \left(\frac{1}{\tilde r} - \tilde r \right) \cos \theta,
	\end{cases}
	\end{equation}
where $\theta$ is the tangent angle at the point $(\tilde z, \tilde r)$. 

The metric is degenerate, so generic geodesics will bounce off as they get close to the $\tilde z$-axis. To obtain smooth embedded caps for the construction, we have to select only the geodesics that tend to the $\tilde z$-axis (and which will eventually become perpendicular to the $\tilde z$-axis). They form a one parameter family of solutions to \eqref{eq:system-ode} characterized by the initial conditions $\tilde z_0= \tilde c$, $\tilde r = 0$, and $\theta_0 = \pi/2$. Because of the metric, the existence and uniqueness of such solutions do not follow from standard ODE methods but from the (un)stable manifold theorem. The flexibility here comes from this one parameter family of rotationally symmetric self-shrinking caps; and a prescribed  unbalancing dictates which cap to select and the radius $\tilde R$ of the intersection circle. 

The asymptotic behavior of our self-shrinkers is also different from the previous constructions in \cite{kapouleas;embedded-minimal-surfaces} and \cite{mine;self-trans}. In both of these articles, the constructed examples tended exponentially fast to the asymptotic catenoids or grim reapers. In this case, although the initial configurations all involve the $\tilde x \tilde y$-plane, the constructed self-shrinkers  are asymptotic to cones at infinity \cite{mine;part2}. 

\subsection{How this construction is similar to previous ones.}

In \cite{kapouleas;embedded-minimal-surfaces} and \cite{mine;self-trans}, the desingularizing surfaces were not only unbalanced but their wings were bent as well to ensure that the solutions to the linearized equation $\lin v =E $ could be adjusted to have exponential decay. The decay is crucial to control the error generated from the cut-off functions when patching up the local solutions  to the linearized equation to form a global solution. In this construction, we can impose an added invariance with respect to the half-turn rotation about the $\tilde x$-axis, and this extra symmetry forces exponential decay on the solutions along the wings of the desingularizing Scherk surface. The situation is similar to the one in \cite{mine;tridents} and we do not need any bending of the wings. 
 
All the estimates and results about the linearized equation on the desingularizing surface $\Sigma$ are obtained by arguments analogous to the ones in \cite{kapouleas;embedded-minimal-surfaces}, although our construction is simpler because there is no bending. Indeed, the difference between equation \eqref{eq:self-shrinker-large} and $H=0$  is of order $\tau$ and the respective linear operators also differ by terms of order at least $\tau$. Since the proofs are very technical and not enlightening, we will not repeat them in this article but just state the relevant properties. The reader who wishes  more details can find some in \cite{mine;tridents}, where we adapted all of the proofs for the equation  $H - \tau \vec e_y \cdot \nu=0$. At this point, we would like to warn the reader that this article is not self-contained and we rely on the reader's familiarity with similar constructions, especially \cite{kapouleas;embedded-minimal-surfaces} or \cite{mine;tridents}, for the proofs of Propositions \ref{prop:4-20} and \ref{prop:7-1}.

Once we define the correct Banach spaces of functions and norms to consider, the few last steps in this article are similar to the ones from Kapouleas's. Namely, the proof that the linearized equation  on the initial surface  $\tilde M$ can be solved modulo a multiple of a well-chosen function follows the same lines as in Kapouleas' article. The final fixed point argument is also similar. Because it would have been strange to stop the construction right before its conclusion, we have included these proofs for the sake of completeness.

\subsection*{Acknowledgments} The author would like to thank Sigurd Angenent for introducing her to this problem and for his encouragement to persist in solving it. 

The author would also like to thank the referees for helpful comments in clarifying the presentation of this article. They have also pointed out that even though all the essential ideas and estimates were present in Section \ref{ssec:linear-P}, the weighted H\"older norms defined in the original version did not yield the compact embeddings necessary to apply the Schauder Fixed Point Theorem. The error has been corrected in this version with the addition of Definition \ref{def:C-cone} and Lemma \ref{lem:cone-op}. 

The original version of this article contained an extra parameter $\varphi$ for bending the wings of the desingularizing surfaces. The imposed invariance under rotation of $180^{\circ}$ around the $\tilde x$-axis makes $\varphi$ superfluous and the author has removed it. 

After completion of the original manuscript, the author learned that N. Kapouleas, S. Kleene, and N. M{\o}ller have announced a similar result \cite{kapouleas-kleene-moller}. They tackled both points mentioned in the previous two paragraphs correctly in their first version. 

\subsection*{Remark}
The notation $\tilde{\mathbf H} = \tilde H \nu$ (which makes $\nu$ the unit inward normal vector for convex surfaces) and the particular scale  (in which the sphere of radius $\sqrt 2$ in $E^3$ is a solution to the self-shrinker equation) follow the conventions of the previous two installments \cite{mine;part1} and \cite{mine;part2}. The scale differs from the scale in Angenent \cite{angenent;doughnuts}, where the sphere of radius $2$ is self-shrinking. It is also worth noting that the orientation of our normal vector is opposite from the one chosen by   Huisken in \cite{huisken;asymptotic-behavior}.

\subsection{Notation} 
\label{ssec:notations}
	\begin{itemize}
	\item $E^3$ is the Euclidean three space equipped with the usual metric.
	\item $\vec e_x$, $\vec e_y$ and $\vec e_z$ are the three coordinate vectors of $E^3$.
	\item We fix once and for all a smooth cut-off function $\psi$, which is increasing, vanishes on $(-\infty, 1/3)$, and is equal to $1$ on $(2/3, \infty)$. We define the functions $\psi[a,b]: \R \to [0,1]$ which transition from $0$ at $a$ to $1$ at $b$ by 
	\[
	\psi[a,b] (s) = \psi \left(\frac{s-a}{b-a}\right).
	\]
	
	\item We often have a function $s$ defined on the surfaces with values in $\R \cup \{ \infty\}$. If $V$ is a subset of such a surface, we use the notation 
	\begin{equation*}
	V_{\leq a} := \{ p \in V: s(p) \leq a \}, \quad V_{\geq a}:= \{ p \in V :s(p) \geq a\}.
	\end{equation*}
	\item $\nu_{S}, g_S, A_S$, and $H_S$ denote, respectively, the oriented unit normal vector, the induced metric, the second fundamental form, and the mean curvature of an immersed surface $S$ in the Euclidean space $E^3$. 
	\item Given a surface $S$ in $E^3$, which is immersed by $X: S \to E^3$ and a $C^1$ function ${\sigma}: S \to \R$, we call the graph of ${\sigma}$ over $S$ the surface given by the immersion $X+{\sigma} \nu$, and denote it by $S_{\sigma}$. We often use $X+{\sigma} \nu$ and its inverse to define projections from $S$ to $S_{\sigma}$, or from $S_{\sigma}$ to $S$, respectively. When we refer to projections from $S$ to $S_{\sigma}$ or from $S_{\sigma}$ to $S$, we always mean these projections. 
	\item Throughout this article, a surface with a tilde $\tilde{S}$ is a surface in the ``smaller" scale, whereas a surface without a tilde denotes its ``larger" version $S = \frac{1}{\tau} \tilde S=\{ (x,y,z) \in E^3 \mid (\tau x, \tau y, \tau z) \in \tilde S\}$, where $\tau$ is a small positive constant. We also use these conventions for geometric quantities, for example, $H$ is the mean curvature of $S$ and $\tilde H$ is the mean curvature of $\tilde S$. However, these notations apply only  loosely to coordinates: we generally use $x,y,z$ when we are working in a ``larger" scale and $\tilde x,\tilde y,\tilde z$ for objects in a ``smaller" scale but these sets of coordinates are not necessarily proportional by a ratio of $\tau$. 
	 \item	We work with the following weighted H\"older norms:
	\begin{equation*}
	\| \phi: C^{k, \alpha}(\Omega, g, f)\| := \sup_{x\in \Omega} f^{-1}(x) \| \phi:C^{k, \alpha}(\Omega\cap B(x), g)\|,
	\end{equation*}
where $\Omega$ is a domain, $g$ is the metric with respect to which we take the $C^{k, \alpha}$ norm, $f$ is the weight function, and $B(x)$ is the geodesic ball centered at $x$ of radius $1$.	
	\end{itemize}
	

\section{Construction of the desingularizing surfaces}
\label{sec:desingularizing-surfaces}


We introduce the Scherk minimal surface and describe how to unbalance, wrap, and bend it to obtain a suitable desingularizing surface. The parameter $\tau$ is a small positive constant which characterizes how much the desingularizing surface will be scaled to fit in the neighborhood of the intersection circle. Although we do not scale the desingularizing surface in this section, $\tau$ still plays a role here as it determines the radius of the circle around which the Scherk surface is wrapped, and how far we truncate our surface. 

\subsection{The Scherk surface}
\label{ssec:scherk}
The Scherk minimal surface $\Sigma_0$  is given by the equation
	\begin{equation*}
	\Sigma_0=\{ (x, y, z) \in E^3 \mid \sin y = \sinh x \sinh z\}.
	\end{equation*}
This surface was discovered by Scherk and is the most symmetric of a one parameter family of minimal surfaces (see \cite{dierkes-hildebrandt;minimal-surfaces}, \cite{karcher;lecture-notes}, or \cite{karcher;minimal-surfaces}). As $x$ ($z$) goes to infinity, $\Sigma_0$ tends exponentially  to the $xy$-plane ($yz$-plane resp.). More precisely, if  we denote by $H^+$ the closed half-plane $H^+ = \{ (s,y) \in \R^2 \mid s \geq 0\}$, we have the following properties. 

\begin{lemma}[Proposition 2.4  \cite{kapouleas;embedded-minimal-surfaces}]
\label{lem:scherk}
For given $\eps \in (0, 10^{-3})$, there are a constant $a = a(\eps)>0$ and  smooth functions $\sigma: H^+ \to \R$ and $F:H^+ \to E^3$ with the following properties:
	\begin{enumerate}
	\item $F(s, y) = (\sigma(s, y), y, s+a) \in \Sigma_0$, 
	\item $\| \sigma : C^5(H^+, g_{H^+}, e^{-s}) \| \leq \eps$,
	\item $\Sigma_0$ is invariant under rotation of $180^\circ$ around the $x$-axis, and reflections across the planes $y =  \frac{\pi}{2} + k\pi$, $k \in \mathbf Z$. 
	\end{enumerate}
\end{lemma}

The constant $\eps$ is a small constant chosen at this point so $a$ is also fixed. 
We call the surface $F(H^+)$ the \emph{top} wing of $\Sigma_0$, and its image under rotation of $180^\circ$ around the $y$-axis is the \emph{bottom} wing. The \emph{outer} wing is the set of points  $\{ (s +a, y , \sigma(s,y))\}$ and the \emph{inner} wing is $\{ (-s -a, y , \sigma(s,y))\}$. We take as standard coordinates the coordinates $(s,y)$ on each of the wings. If a point of $\Sigma_0$ does not belong to any of the wings, we take its $s$-coordinate to be zero.

\subsection{Unbalancing}


Because the equation $H+\tau^2 X\cdot \nu=0$ is a perturbation of $H=0$, one expects that the respective linear operators $\lin v = \Delta v + |A|^2 v - \tau^2 X \cdot \nabla v + \tau^2 v$ and $Lv := \Delta v + |A|^2$ have similar properties.  The mean curvature is invariant under translations, therefore the functions $\vec e_x \cdot \nu$, $\vec e_y \cdot \nu$, and $\vec e_z \cdot \nu$ are in the kernel of the linear operator $L$ associated to normal perturbations of $H$.  We can rule out $\vec e_y \cdot \nu$ and $\vec e_z \cdot \nu$ by imposing symmetries (see (iii) of Lemma \ref{lem:scherk}). The remaining function $\vec e_x \cdot \nu$ does not have the required exponential decay, however, it indicates  that $L$ has an approximate kernel generated by a function close to $\vec e_x \cdot \nu$ and one can only solve the equation $L v =E $ with a reasonable estimate on $v$ if $E$ is perpendicular to the approximate kernel. But we do not have such control over the inhomogeneous term, so we have to introduce a function $w$ to cancel any component  parallel to the approximate kernel. The function $w$ has to be in the direction of $\vec e_x \cdot \nu$ in the sense that $\int w (\vec e_x \cdot \nu) \neq 0$.

Let $S$ be one period of the desingularizing surface $\Sigma$. According to the balancing formula from \cite{korevaar-kusner-solomon}, the mean curvature of $\Sigma$ satisfies
	\[
	\int_{S}  H  \vec e_x \cdot \nu dg_{S} = 2 \pi \sum_{i=1}^4 v_i \cdot \vec e_x,
	\]
where $v_i$ is the direction of the  plane asymptotic to the $i$th wing. The idea is to define $w$  as a derivative of $H$ and use unbalancing to move  the top and bottom wings toward $\vec e_x$ to generate a multiple of $w$.

\begin{definition}
\label{def:Zb}
For $b \in (-\frac{1}{10}, \frac{1}{10})$, we take a family of diffeomorphisms $Z_{b}: E^3 \to E^3$ depending smoothly on $b$ and satisfying the following conditions:
	\begin{enumerate}
	\item $Z_{b}$ is the identity in the region $\{(x, y, z)\mid |z| < \frac{1}{2}|x| \}$,
	\item in the region $\{ (x, y, z) \mid |z| > 2 |x| \text{ and } |z| > 1\}$,  $Z_{b}$ is the rotation by an angle $b$ about the $y$-axis toward the positive $x$-axis,
	\item $Z_0$ is the identity.
	\end{enumerate}
\end{definition}	

We denote by $\Sigma_b$ the unbalanced surface $Z_b(\Sigma_0)$ and push forward the coordinates $(s,y)$ of $\Sigma_0$ onto this new surface using $Z_b$.  


\subsection{Wrapping the Scherk surface around a circle}
Given $\tilde R \in (1,2)$, let $R = \tau^{-1} \tilde R$ and define the maps $\Phi_{R}:E^3 \to E^3$ by 
	\[
	\Phi_{R}(x, y, z)=R 
			\left(e^{ x/R} \cos(\tau y/\sqrt 2) , e^{x/R} \sin(\tau y/\sqrt 2), z/R\right).
	\]
We allow $\tilde R$ to differ slightly from $\sqrt 2$ so that it can be chosen to fit the self-shrinking rotationally symmetric caps in Section \ref{ssec:fitting-caps}. The exponential factor $e^{x/R}$ has no real significance here and can be replaced by any strictly increasing diffeomorphism $f:\R \to (0,\infty)$ with the property $f(0)=1$. 

We push forward the coordinates $(s, y)$ of $\Sigma_0$ onto the surface $ \Phi_{R} \circ Z_b (\Sigma_0)$. The piece of surface corresponding to $\{ s \leq 0\}$ within the slab $\{ -a \leq z \leq a\}$ is called the \emph{core} of the desingularizing surface and will no longer be modified. The image of the plane asymptotic to the top wing of $\Sigma_b$ is given by 
	\begin{equation}
	\label{eq:asymp-plane}
	\left \{ \left (R e^{z \tan b/ R}\cos ( \tau y/\sqrt 2), R e^{ z \tan b/ R}\sin ( \tau y/\sqrt 2) , z \right), z \geq a \right\},
	\end{equation}
where we placed the boundary of the asymptotic surface at $z=a$ to simplify subsequent computations. In general, the piece of sphere given by \eqref{eq:asymp-plane} does not satisfy the equation for self-shrinkers \eqref{eq:self-shrinker} and would generate too much error if it were used to build the top wing. Instead, we record its boundary and conormal direction in Definition \ref{def:pivot} and equation \eqref{eq:beta}, respectively, then choose a piece of rotationally self-shrinking surface these boundary conditions in Definition \ref{def:kappa}.


 \begin{definition} 
 \label{def:pivot}
The circle that bounds the asymptotic surface given by \eqref{eq:asymp-plane} is called the \emph{pivot} of the top wing.
\end{definition}
We define $\beta $ to be the angle the inward conormal vector makes with the direction $\vec e_z$ at the pivot, which is given by the equation
	\begin{equation}
	\label{eq:beta}
	\tan \beta = (\tan b) e^{ a \tan b/R}.
	\end{equation}
Note that $\beta$ is a smooth function of $b$ and that
	$
	1-C\tau \leq \left|\frac{d\beta}{db}\right| \leq 1+C\tau
	$	
for some positive constant $C$.

\begin{definition}
\label{def:kappa}
Given $\tau$ a small positive constant, $b \in (-\frac{1}{10}, \frac{1}{10})$ and $\tilde R \in [1.3, 1.5]$, we consider the solution $(\tilde z(t), \tilde r(t), \theta(t))$ to the system \eqref{eq:system-ode} with initial conditions 
	\[
	\tilde z(0) = \tau a, \quad \tilde r(0)= \tilde R e^{a \tau \tan b/\tilde R}, \quad  \theta(0) =\beta,
	\]
where $\beta$ is given by \eqref{eq:beta}.
We define the map $\kappa[\tilde R, b, \tau]: H^+ \to E^3$ by 
	\[
	\kappa[\tilde R, b, \tau](s, y) = \frac{1}{\tau}\left(\tilde r (t(s)) \cos (\tau y/\sqrt 2), \tilde r (t(s)) \sin(\tau y/\sqrt 2), \tilde z(t(s)) \right),
	\]
where we reparamatrize using $t(s)$ satisfying $\frac{dt}{ds} = \tau \tilde r(t)/\sqrt 2$, $t(0) =0$ so that $\kappa[\tilde R, b, \tau]$ is conformal. 
\end{definition}

From standard ODE theory, the flow of \eqref{eq:system-ode} is smooth and depends smoothly on  the initial conditions. The surfaces $\kappa[\tilde R,b, \tau](H^+)$ are not  embedded in general, but we only consider the small pieces where $0 \leq s \leq 5 \delta_s/\tau$, where the small positive constant $\delta_s$ will be determined in Section \ref{sec:estimates}. The pull-back of the induced metric by $\kappa$ is $\rho^2(ds^2+dy^2)$ where $\rho^2= \tilde r^2 /2.$

\subsection{The inner and outer wings} The construction of these two wings is very simple: we just use the transition function $\psi[4 \delta_s/\tau, 3 \delta_s/\tau] \circ s$ to cut off the graph of  $\sigma$ over these two wings, then truncate the desingularizing surface at $s=5 \delta_s/\tau$.


\subsection{The top and bottom wings.}

The top wing of our desingularizing surface is essentially the graph of $\sigma$ from Lemma \ref{lem:scherk} over the self-shrinking rotationally symmetric surface $\kappa[\tilde R, b, \tau] (H^+)$. To get smooth transitions, we use the cut-off functions $\psi[1,0]$ and $\psi_s$. The factor $\tau^{-1}$ in the scale of $\psi_s$ is to ensure that we cut off $\sigma$ far enough as not to generate too much error, while $\delta_s$ is small in order to control the geometry of the wings and so that the desingularizing surfaces fit in a small neighborhood of the intersection circle (in the original scale of \eqref{eq:self-shrinker}). 
	\begin{definition}
	\label{def:F-wing}
	For given $\tau, \tilde R$, and $b$ as in Definition \ref{def:kappa}, we define $F[\tilde R,b,\tau]: H^+ \to E^3$ by 
	\begin{multline*}
	F[\tilde R, b, \tau] (s,y)= \psi[1,0](s)  \Phi_{\tilde R/\tau} \circ Z_b \circ F (s, y) \\
				+ \big(1 - \psi[1,0](s)\big) \big( \kappa[\tilde R, b,\tau](s,y) + \psi_s(s)\sigma(s,y) \nu[\tilde R, b, \tau](s,y) \big),
	\end{multline*}
where $\psi_s$ is defined by $\psi_s(s) = \psi[4 \delta_s /\tau, 3 \delta_s /\tau](s)$ and $\nu[\tilde R, b, \tau]$ is the Gauss map of $\kappa[\tilde R, b, \tau](H^+)$ chosen so that $\nu[\tilde R, b,\tau] (0,0) = (\cos(\beta), 0 , -\sin(\beta))$
	\end{definition}
The top wing is divided into four regions: 
	\begin{itemize}
	\item $\{ 0 \leq s \leq 1\}$ is a transition region from the core to the bent wing.
	\item on $\{ 1 \leq s \leq 3\delta_s/\tau\}$, the wing is the graph of $\sigma$ over the asymptotic rotationally symmetric piece of self-shrinker.
	\item $\{ 3 \delta_s/\tau \leq s \leq 4 \delta_s/\tau\}$ is another transition region where we cut off the graph of $\sigma$.
	\item on $\{ 4 \delta_s /\tau \leq s \}$, the wing is a piece of rotationally symmetric self-shrinker.
	\end{itemize} 
As before, we push forward the coordinates $(s,y)$ of $H^+$ onto  $F[\tilde R, b, \tau](H^+)$ and clip the wing at $s=5\delta_s/\tau$. The bottom wing is obtained by rotating the top wing by $180^{\circ}$  around the $x$-axis. 
\subsection{The desingularizing surfaces}
The desingularizing surface $\Sigma[\tilde R, b, \tau]$ is composed of the core, the inner and outer wings, and the top and bottom wings defined in the three previous sections. We sometimes denote it by  $\Sigma$ for simplicity. The next proposition collects some useful properties of our desingularizing surfaces. 
	
\begin{proposition}
\label{prop:smooth-dependence}
There exists a constant $\delta_{\tau}>0$ such that for $\tau \in (0,\delta_{\tau})$, $\tilde R \in [1.3,1.5]$, and $b \in (-\frac{1}{10}, \frac{1}{10})$, the surface $\Sigma[\tilde R, b,  \tau]$ satisfies the following properties:
	\begin{enumerate}
	\item $\Sigma[\tilde R, b,  \tau]$ is a smooth surface immersed in $E^3$ which depends smoothly on its parameters.
	\item If $\tau=\frac{\sqrt 2}{m}$, $m \in \mathbf N$, the surface $\Sigma[\tilde R, b, \tau]$ is embedded. Moreover, $\Sigma[\tilde R, b,\tau]$ is invariant under the rotation of $180^{\circ}$ around the $x$-axis and under the reflections across planes containing the $z$-axis and forming angles $\frac{\pi}{2m} + \frac {k \pi}{m}$,  $k \in \mathbf Z$, with the $xz$-plane. 
	\end{enumerate}
\end{proposition}

We choose a positive integer $m$ so that $\tau = \sqrt 2/m \in (0, \delta_{\tau})$ and fix the value of $\tau$ for Sections \ref{sec:estimates}, \ref{sec:construction}, \ref{sec:linear}, and \ref{sec:quadratic}.


\section{Estimates on the desingularizing surfaces}
\label{sec:estimates}

In this section, we claim that the desingularizing surfaces $\Sigma$ are suitable approximate solutions. All the estimates from Section 4 in \cite{kapouleas;embedded-minimal-surfaces}  are valid, with $H$ replaced by $H_{\Sigma}+\tau^2 X_{\Sigma} \cdot \nu_{\Sigma}$ and the corresponding linear operator $L_S = \Delta_S + |A_S|^2$ replaced by $\lin_S= \Delta_S + |A_S|^2 + \tau^2 (1-X\cdot \nabla)$. The factor $\tau^2$ combats the scale of the position $X \sim \tau^{-1}$ so the extra term does not add significantly. This is where we choose the constant $\delta_s$ small enough so that the metrics $g_{\Sigma}$, $\kappa_{\ast} (ds^2 + dy^2)$, and $F[\tilde R, b, \tau]_{\ast}(ds^2+dy^2)$  are uniformly equivalent on the top and bottom wings. Later on, in the proof of Proposition \ref{prop:7-1}, $\delta_s$ may be changed to an even smaller constant so that the linear operator $\lin_{\Sigma}$ on the wings can be treated as a perturbation of the Laplace operator on a flat cylinder. Because the proofs are technical and do not showcase the main aspects of the construction, we will not reproduce them here. The reader can find all the details in Section 4 of \cite{kapouleas;embedded-minimal-surfaces} or Section 3 of \cite{mine;tridents}.

In what follows, the parameter $\tau$ and the radius $\tilde R$ are fixed and the dependence on $\tilde R$ will be omitted. Moreover, because $\tilde R$ takes values in a compact set, all of the constants $C$ can be chosen independently of $\tilde R$.

We define the function $w: \Sigma_0 \to \R$ by 
	\[
	w:= \left. \frac{d}{db} \right|_{b=0} H_b \circ Z_b,
	\]
where $H_b$ denotes the mean curvature on the surface $Z_b(\Sigma_0)$. 

The main contribution to $H_{\Sigma} + \tau^2 X_{\Sigma} \cdot \nu_{\Sigma}$ comes from the unbalancing term ($bw$). Here $\gamma$ is a constant in $(0,1)$ which indicates that the exponential decay is slower due to the presence of the cut-off function $\psi_s$ in Definition \ref{def:F-wing}.

\begin{proposition}
\label{prop:4-20}

For $(b, \tau)$ as in Proposition \ref{prop:smooth-dependence}, the quantity $H_{\Sigma}+\tau^2 X_{\Sigma} \cdot \nu_{\Sigma}$ on $\Sigma = \Sigma[\tilde R, b, \tau]$ satisfies 
	\[
	\| H_{\Sigma}+\tau^2 X_{\Sigma} \cdot \nu_{\Sigma} - b w :C^{0, \alpha}(\Sigma, g_{\Sigma}, e^{-\gamma s}) \| \leq C(\tau + |b|^2). 
	\]
\end{proposition}


\section{Construction of the initial surfaces}
\label{sec:construction}

In the construction of the desingularizing surfaces, we did not unbalance or bend the inner and outer wings so attaching them to a disk and plane respectively is straightforward. For the top and bottom wings, the story is more complicated. In the case of minimal surfaces, coaxial catenoids form a two parameter family of minimal surfaces whose embeddings depend smoothly on the parameters. When the desired tangent direction of a gluing wing is changed, one has the flexibility of attaching a catenoid close to the original one. To get flexibility in this construction, we have to consider the sphere of radius $\sqrt 2$ as a member of a family of self-shrinking surfaces and not as a member of a family of spheres. Note that in \cite{kapouleas;embedded-minimal-surfaces}, Kapouleas had invariance for reflection across planes only, so he used a two-parameter family of initial configurations. Since we have an additional symmetry (invariance under rotation of $180^{\circ}$ around the $x$-axis), the family of self-shrinking surfaces depends  on one parameter only. 

\subsection{A family of rotationally symmetric self-shrinking caps} 
\label{ssec:shrinking-caps}

 In \cite{angenent;doughnuts}, Angenent showed that hypersurfaces of revolution are self-shrinkers if and only if they are generated by  geodesics of the half-plane $\{ (\tilde z, \tilde r) \mid \tilde r \geq 0\}$ equipped with the metric \eqref{eq:metric}. Given any point $(\tilde z, \tilde r)$ and an angle $\theta \in \R$, there is a unique geodesic  through $(\tilde z, \tilde r)$ with tangent vector $(\cos \theta, \sin \theta)$. Such a geodesic parametrized by arc length satisfies the following system of ODEs
	\begin{equation}
	\tag{\ref{eq:system-ode}}
	\begin{cases}
	\dot {\tilde z} = \cos \theta, \\
	\dot {\tilde r} = \sin \theta, \\
	\dot \theta = \tilde z \sin \theta + \left(\frac{1}{\tilde r} - \tilde r \right) \cos \theta.
	\end{cases}
	\end{equation}

 Because the metric becomes degenerate as $\tilde r \to 0$, geodesics in general will bounce off as they approach the $\tilde z$-axis. For the purpose of having a complete rotationally symmetric cap, we will only consider geodesics that tend towards  the $\tilde z$-axis. Such curves will always meet the $\tilde z$-axis at a right angle. 
\begin{definition}
\label{def:ode-flow}
For $ c \in \R$ close to $\sqrt 2$, we denote by $\tilde \gamma_{ c}(\cdot)$ or $\tilde \gamma( c; \cdot)$ the geodesic in the half-plane $\{ (\tilde z, \tilde r) \mid \tilde r \geq 0\}$ equipped with the metric \eqref{eq:metric} with initial conditions 
	\[
	\tilde \gamma( c; 0) = ( c, 0), \quad \tilde \gamma'( c; 0) = (0, 1).
	\]
\end{definition}
Note that the curve $\tilde \gamma_{ c}$ is the projection to the $\tilde z\tilde r$-plane of the solution  to \eqref{eq:system-ode}  $\tilde \alpha_{ c}(\cdot)=(\tilde z( c;\cdot), \tilde r( c;\cdot),   \theta(c;\cdot))$ with initial conditions 
	\[
	\tilde z( c;0)= c, \quad \tilde r( c;0) = 0, \quad  \theta( c;0) = \pi/2.
	\]
The solution corresponding to the hemisphere of radius $\sqrt 2$ is 
	\[
	\tilde z (\sqrt 2; t) = \sqrt 2 \sin\left(\frac{\pi}{2} + \frac{t}{\sqrt 2} \right),  
	\tilde r (\sqrt 2; t) =- \sqrt 2 \cos\left(\frac{\pi}{2} + \frac{t}{\sqrt 2} \right),
	\theta (\sqrt 2;t ) = \frac{\pi}{2}+ \frac{t}{\sqrt 2}.
	\]

\begin{proposition}
\label{prop:smooth-dependence-orbits}
There exists a constant $\delta_c>0$ for which the map $(\tilde z,\tilde r,\theta): (\sqrt 2 - \delta_c, \sqrt 2 + \delta_c) \times [0,  3 \pi/\sqrt 2] \to \R^3$ that associates $( c,t)$ to $(\tilde z( c;t), \tilde r( c;t),  \theta( c;t))$ in Definition \ref{def:ode-flow} is smooth. 
\end{proposition}

The number $3\pi/\sqrt 2$ was chosen so that all the geodesics $\tilde \gamma_{c}$ would exist long enough to exit the first quadrant.

\begin{proof}
The system of ODEs \eqref{eq:system-ode} can be reparametrized using the variable $h$ for which $\frac{d}{dh} = \tilde r \frac{d}{dt}$ to get
	\[
	\begin{cases}
	\frac{d}{dh}\tilde  z = \tilde r \cos \theta, \\
	\frac{d}{dh} \tilde r = \tilde r \sin \theta,\\
	\frac{d}{dh}  \theta = \tilde z \tilde r \sin \theta + (1 - \tilde r^2) \cos  \theta.
	\end{cases}
	\]
With this parametrization in the octant $\{ \tilde r,\tilde z,\theta \geq 0\}$, 
	\begin{itemize}
	\item the line $\{\tilde z=0, 0\leq \tilde r< \infty, \theta = \pi/2\}$ is invariant (and the corresponding self-shrinker is the plane),
	\item the line $\{0\leq \tilde z < \infty, \tilde r=1,  \theta=0\}$ is invariant (and the corresponding self-shrinker is the cylinder),
	\item the line $l =\{0\leq \tilde z < \infty, \tilde r=0,  \theta=\pi/2\}$ consists of fixed points. The linearization at $(c, 0, \pi/2) \in l$ is 
	\[
	\frac{d}{dh} 
	\left(
	\begin{array}{c}
	\delta \tilde z\\
	\delta \tilde r\\
	\delta \theta
	\end{array}
	\right)
	=
	\left(
	\begin{array}{ccc}
	0 & 0 & 0 \\
	0 & 1 & 0 \\
	0 & c & -1 
	\end{array}
	\right) 
	\left(
	\begin{array}{c}
	\delta \tilde z\\
	\delta \tilde r\\
	\delta  \theta
	\end{array}
	\right).
	\]
	\end{itemize}
The line $l$ is therefore normally hyperbolic, with a stable manifold contained in the plane $\{ \tilde r=0\}$. Its unstable manifold consists of a one parameter family of orbits $\tilde \alpha_c$, each $\tilde \alpha_c$ emanating from the point $(  c, 0, \pi/2)$. The short-time existence and uniqueness of these orbits $\tilde \alpha_c$, as well as the smoothness of the unstable manifold is given by the (un)stable manifold theorem (see  Theorem (4.1) \cite{hirsch-pugh-shub} or  Theorem III.8 \cite{shub;global-stability}). Once we get away from the line $l$ using this first step, we can  extend the one parameter family of orbits $\tilde \alpha_c$ smoothly using standard ODE theory. The uniform dependence of $(\tilde z, \tilde r, \theta)$ for $t \in [0, 3 \pi/\sqrt 2]$ is obtained by a compactness argument since the system \eqref{eq:system-ode} is not singular for $\tilde r$ away from zero. 
\end{proof}

\begin{proposition} 
\label{prop:dependence-theta-c}
There exists a positive constant $\delta_{\theta}$ such that given $\theta_0 \in (\pi - \delta_{\theta}, \pi+\delta_{\theta})$ there exists a unique constant $ c_0 \in (\sqrt 2-\delta_c, \sqrt 2+\delta_c)$ for which the orbit $\tilde \alpha_{ c_0}$  hits the $\tilde r$-axis at an angle $\theta_0$. Moreover, for some constant $C$ independent of $ \theta_0$, we have
		\[
		|  c_0 - \sqrt 2| < C | \theta_0 - \pi|.
		\]
\end{proposition}

\begin{proof} We start the proof by giving a different description of the geodesics. The graph of a function $f$ over the circle of radius $\sqrt 2$, i.e. the curve given by $(\sqrt 2 + f(t)) (\cos t,  \sin t )$, is a self-shrinker if and only if
		\[
		\frac{-f'^2+f''(\sqrt 2 + f)}{f'^2 +(\sqrt 2 + f)^2 } + \frac{f' \cos t}{(\sqrt 2+f) \sin t} + (\sqrt 2 +f)^2 -2 =0.
		\]
Using the change of variable $h(t) = \ln (\sqrt 2 +f(t))$, the equation above is equivalent to
		\begin{equation*}
		\frac{h''}{1+h'^2} + \frac{\cos t}{\sin t} h'+ e^{2h}-2=0.
		\end{equation*}

The existence of a solution $h( c; t)$ with $h(0)= c$, $h'(0)=0$ follows from the proof of Proposition \ref{prop:smooth-dependence-orbits}. In addition, the unstable manifold theorem gives the smooth dependence of the solution $h$ on its parameters and  
	\[
	h( c; t) = \ln(\sqrt 2)+ \psi(t) ( c-\sqrt 2) + o( c-\sqrt 2),
	\]
where $o(\eps) \to 0$ as $\eps \to 0$ and where the function $\psi(t)$ satisfies the linear ODE
	\[
	\psi''+ \frac{\cos t}{\sin t} \psi'+ 4 \psi =0, \quad \psi(0)=1, \psi'(0) = 0.
	\]
The solution is given by $\psi(t)= P_{\frac{1}{2}(-1+\sqrt 17)} (\cos t)$, where $P_{\lambda}(t)$ is the Legendre function. We get the existence of $c_0$ and the estimate \eqref{eq:dependence-theta-c} below because the derivative $dP/dt$ is positive at $t=\pi/2$.
\end{proof}

In the following corollary, we seek to hit the line $\tilde z = \tau a$ at a specific angle $\theta_1$. 

\begin{corollary}
\label{cor:dependence-theta-c}
There exists a positive constant $\delta_{\theta}$ independent of $\tau$ such that given $\theta_1$ with 
	\[
	| \theta_1 - (\pi-\sin^{-1}(\tau a/\sqrt 2))| \leq \delta_{\theta},
	\]
 there exists a unique constant $ c_1\in (\sqrt 2 - \delta_c, \sqrt 2+\delta_c)$ for which the orbit $\tilde \alpha_{ c_1}$ hits the line $\tilde z=\tau a $ at an angle $\theta_1$ and 
	\begin{equation}
	\label{eq:dependence-theta-c}
	| c_1 - \sqrt 2| \leq C | \theta_1 - \pi + \sin^{-1}(\tau a/\sqrt 2)|
	\end{equation}
\end{corollary}
\begin{proof}
Because $\tau$ is a small constant, this corollary follows from Proposition \ref{prop:dependence-theta-c} and the smooth dependence on $c$ from Proposition \ref{prop:smooth-dependence-orbits}.
\end{proof}

\subsection{Fitting the self-shrinking caps to the desingularizing surfaces}
\label{ssec:fitting-caps}

Let us recall that in Section \ref{sec:desingularizing-surfaces}, we did not restrict ourselves to geodesics that meet the $\tilde z$-axis perpendicularly but considered any solution to \eqref{eq:system-ode} to construct the asymptotic surfaces $\kappa[\tilde R, b, \tau](H^+)$. We now choose the radius $\tilde R$  in function of the angle $b$ so that the surface asymptotic to the top wing of $\Sigma$ is contained in a self-shrinking rotationally symmetric cap. Given $b$, we take $\tilde R(b)$ to be the $c_1$ given in Corollary \ref{cor:dependence-theta-c} corresponding to $\theta_1 = \beta$, where $\beta$ is given by \eqref{eq:beta}.

\subsection{Construction of the initial surfaces $\tilde M(b, \tau)$}
\label{ssec:initial-surface}

Let us recall that $\tau = \sqrt 2/m \in (0, \delta_{\tau})$, for a previously chosen integer $m$. We fix a constant $\zeta$ which will be determined in the proof of Theorem \ref{thm:8-2}.

Given $b \in [-\zeta \tau, \zeta \tau]$, we start the construction of the initial surface by taking the desingularizing surface $\Sigma[\tilde R(b), b, \tau]$ and shrinking it to $\tilde \Sigma=\tilde \Sigma[\tilde R(b), b, \tau]$ with the homothety $\hcal$ of ratio $\tau$ centered at the origin. We top off (on the top and bottom) the desingularizing surface $\tilde \Sigma$ with self-shrinking caps generated by rotating the curve $\tilde \gamma_{ c(b)}$ around the $\tilde z$-axis. The inner wing of $\tilde \Sigma$ is attached to a flat disk and the outer wing to a plane. 

\begin{definition}
\label{def:components}
The surface constructed in the above paragraph is denoted by $\tilde M(b, \tau)$. We push forward the function $s$ by $\hcal$ from $\Sigma$ to $\tilde \Sigma$ and extend it to the whole surface $\tilde M(b, \tau)$ by taking $s=5\delta_s/\tau$ on $\tilde M\setminus \tilde \Sigma$. 

Let $\ubar a := 8 |\log \tau|$. We define 
	\begin{gather*}
	\tilde \dcal=\textrm {the component of $\tilde M_{\geq \ubar a}$ that contains the inner disk}\\
	\tilde \pcal=\textrm {the component of $\tilde M_{\geq \ubar a}$ that contains the outer plane}\\
	\tilde \ccal= \tilde M_{\geq \ubar a} \setminus (\tilde \dcal \cup \tilde \pcal)
	\end{gather*}
and their image under $\hcal^{-1}$ by $\dcal, \pcal$, and $\ccal$ respectively.
\end{definition}
Choosing the constant $\ubar a$ of order $|\log \tau|$ will be useful for getting a contraction (equation \eqref{eq:contraction}) in the proof of Theorem \ref{thm:linear-invertible}.
\begin{proposition}
\label{prop:smooth-initial-surfaces}
Given a positive integer $m$ so that $\tau = \sqrt 2/m \in (0, \delta_{\tau}) $ and $b \in [-\zeta \tau, \zeta \tau]$, the surface $\tilde M=\tilde M(b, \tau)$ is well defined by the construction above and satisfies the following properties:
\begin{enumerate}
\item $\tilde M$ is a complete smooth embedded surface which depends smoothly on $(b)$.
\item $\tilde M$ is invariant under rotation of $180^{\circ}$ around the $\tilde x$-axis.
\item $\tilde M$ is invariant under the action of the group $G$ generated by reflections across the planes containing the $\tilde z$ axis and forming an angle of $\frac{ \pi}{2m}+k \frac{\pi}{m}$, $k \in \mathbf Z$ with the $\tilde x \tilde z$-plane.
\item As $m \to \infty$, the sequence of initial surfaces $\tilde M(b, \tau)$ converges uniformly in $C^k$, for all $k \in \mathbf N$, to the union of a sphere of radius $\sqrt 2$ and the $\tilde x \tilde y$-plane on any compact subset of the complement of the intersection circle.
\item Let us denote by $T$ the translation by the vector $-\sqrt 2 \vec e_x$. As $m \to \infty$, the sequence of surfaces $m T(\tilde M(b,  \tau))$ converges uniformly in $C^k$, for all $k \in \mathbf N$, to the Scherk surface $ \Sigma_0$ on any compact subset of $E^3$. 
\end{enumerate}

\end{proposition}

The parameter $\tau$ will always be equal to $\sqrt 2/m$ for some natural number $m$ from now on. 

\section{The linearized equation}
\label{sec:linear}

We study the linearized equations on the various pieces $\Sigma$, $\tilde \ccal$, $\tilde \dcal$, and $\tilde \pcal$ and find appropriate estimates for the solutions. The linearized equation on the whole surface $M$ is solved by using cut-off functions to restrict ourselves to the various pieces and patching up all these local solutions with an iteration process.  

\subsection{The linearized equation on $\Sigma$}

The linear equation $\lin_{\Sigma} v = \Delta_{\Sigma} v +|A_{\Sigma}|^2 v + \tau^2 v - \tau ^2 X \cdot  \nabla v = E$ on $\Sigma=\Sigma[\tilde R, b, \tau]$ can be solved modulo the addition of a multiple of $w$ on the right hand side, which takes care of small eigenvalues of $\lin$. The next proposition is reminiscent of Proposition 7.1 in \cite{kapouleas;embedded-minimal-surfaces} but the proof is simpler. In our case, the group of imposed symmetries is larger and is used to rule out troublesome linear growth for solutions of $\Delta v=E$ (or $\lin v=E$). The exponential decay is therefore achieved without resorting to any bending or adjustment along the wings. Proposition \ref{prop:7-1} below is similar to Corollary 22 in \cite{mine;tridents}. Because one can prove it by following the steps in \cite{mine;tridents} and simply substituting the linear operator, we do not give details of the proof here, only the main ideas.
\begin{proposition}
\label{prop:7-1}
Given $E' \in C^{0,\alpha}(\Sigma)$, there are $b_{E'}  \in \R$ and $v_{E'} \in C^{2,\alpha}(\Sigma)$ such that:
	\begin{enumerate}
	\item $b_{E'}$ and $v_{E'}$ are uniquely determined by the proof, 
	\item $\lin_{\Sigma} v_{E'} = {E'} + b_{E'} w $ on $\Sigma$ and $v_{E'}=0$ on $\pd \Sigma$,
	\item $|b_{E'}| \leq C \| E'\|$, where $\|E' :C^{0,\alpha}(\Sigma, g_{\Sigma}, e^{-\gamma s} ) \|$, 
	\item $ \| v_{E'}: C^{2,\alpha} (\Sigma, g_{\Sigma}, e^{-\gamma s})\| \leq C \| E'\|$.
	\end{enumerate}
\end{proposition}

\begin{proof}[Sketch of the proof] It suffices to prove the result for the operator $L_{0}=\Delta_{\Sigma_0} +|A_{\Sigma_0}|^2$ on the piece of the original Scherk surface $\Sigma_{0, \leq 5\delta_s/\tau}$. Indeed, because $\tau$ and $\delta_s$ are small constants, one period of our desingularizing surface $\Sigma$ is diffeomorphic to $\Sigma_{0, \leq 5\delta_s/\tau}$, their respective metrics are uniformly equivalent, and we can control $|A_{\Sigma}|^2-|A_{\Sigma_0}|^2$. Moreover, the linear operator $\lin_{\Sigma}$ on the wings can be treated as a perturbation of the Laplace operator on long flat cylinders (of length $l=5 \delta_s/\tau$) because $ l^2  |A_{\Sigma}|^2$ is small if $\delta_s$ is small. 

The first step is to show that the linear operator $L_{0}$ on $\Sigma_{0,\leq 5 \delta_s/\tau}$ is Fredholm of index $-1$. Here, when we say  ``the operator $L_0$ on a surface $S \subset \Sigma_0$", we consider $L_0$ as an operator from the space $C_0^{2,\alpha}(S, g_{0}, e^{-\gamma s})$ of exponentially decaying $C^{2,\alpha}$ functions with vanishing boundary conditions on $\pd S$ to $C^{0,\alpha}(S, g_{0}, e^{-\gamma s})$.  To compute the Fredholm index, we use a Mayer-Vietoris-type argument. We split $\Sigma_{0, \leq 5\delta_s/\tau}$ into its four wings and a slightly bigger core, $\Sigma_{0, \leq 2}$, and show that the Fredholm index on the whole surface is the sum of the indices on each piece minus the sum of the indices on the overlaps. 

The overlaps are cylinders of length $2$. Like the Laplace operator, $L_0$ is invertible and has index $0$. On the wings, $L_0$ can also be compared to the Laplacian, but the situation is more complicated because of the exponential decay. A computation with Fourier series shows that the exponentially decaying solution to 
$ \Delta v (s, y)=E(s, y)$, $v|_{s=5 \delta_s/\tau}=0$ vanishes at $s=0$ if and only if 
	\begin{equation}
	\label{eq:ortho}
	\int_0^l \int_0^{2 \pi} s E(s, y) dy ds=0.
	\end{equation}
Because all the surfaces are invariant under rotation of $180^{\circ}$ around the $x$-axis, the condition \eqref{eq:ortho} is satisfied on the inner and outer wings.  Also, functions on the bottom wings are tied to functions on the top wings so the Fredholm index of $L_0$ on the union of the four wings is $-1$. 
 
The Scherk surface is a minimal surface; therefore the Gauss map $\nu: \Sigma_0 \to S^2$ is conformal. It sends $\Sigma_{0,\leq 2}$ to the sphere minus four discs centered at $(\pm 1, 0, 0)$ and $(0, 0, \pm 1)$, and the pull-back $\nu^{\ast}g_{S^2}$ of the metric on $S^2$ is $\nu^{\ast}g_{S^2}=\frac{1}{2}|A_{\Sigma_0}|^2 g_0$. This means that a function on $\Sigma_0$ is in the kernel of $L$ if its pushforward by the Gauss map is in the kernel of $\Delta_{g_{S^2}}+2$ on the sphere. For domains $U\subset S^2$ with smooth boundary, the Laplace operator from $\{ f \in C^{2,\alpha}(U)\mid f|_{\pd U} =0\}$ to $C^{0,\alpha}(U)$ is invertible. The addition of a compact operator does not change the index, so $\Delta + 2$ on $U$ and $L_0$ on $\Sigma_{0, \leq 2}$  both have index $0$.

The operator $L$ on $\Sigma_{0, \leq 5 \delta_s/\tau}$ therefore has index $-1$. Without the additional symmetry, the index would be $-4$. We would need to insert more parameters in the construction of the initial surface to account for the bigger cokernel, hence the bending of the wings. The bending is there to cancel possible linear contributions from the kernel of $\Delta$ on each wing. If the original configuration has the extra half-turn symmetry as in this case or the case a plane and a grim reaper in \cite{mine;tridents}, it is not needed.

Let us recall that $L$ is the linear operator associated to normal perturbations of the mean curvature by graphs of functions. Because the mean curvature is invariant when a surface is translated, the functions $\vec e_x \cdot \nu$, $\vec e_y \cdot \nu$, and $\vec e_z \cdot \nu$ are in the kernel of $L$. We can rule out the last two with the imposed symmetries on $\Sigma_0$, namely that any perturbation of $\Sigma_0$ should be invariant under rotation of $180^{\circ}$ around the $x$-axis and reflections across the planes $y=\frac{\pi}{2}+k \pi$, $k \in \mathbf Z$. The remaining one, $\vec e_x \cdot \nu$, does not decay exponentially, so it is not in the kernel. Nevertheless, it plays a role here in characterizing the cokernel:  the map  $\Phi:\mathcal C^{0, \alpha}(\Sigma_{0, \leq 5 \delta_s/\tau}, g_0, e^{-\gamma s}) \to \R$ defined by 
	\[
	\Phi(h) = \int_{\Sigma_{0, \leq 5 \delta_s/\tau}} h (e_1 \cdot \nu) d\mu,
	\]
where $d\mu$ is the Hausdorff measure on $\Sigma_0$, vanishes on the range of $L$ and  $\Phi(w) \neq 0$. The cokernel therefore has dimension 1 and the kernel is trivial. Given a function $E'$, one can take $b_{E'} = - \Phi(E') /\Phi(w)$ and find a solution $v$ to the Dirichlet problem $L v=E'+b_{E'}w$, $v=0$ on $\pd \Sigma_{0,\leq 5 \delta_s/\tau}$. 
\end{proof}

\subsection{The linearized equation on the cap $\tilde \ccal$}

Because this surface is in the ``smaller" scale, we consider the linear operator $\tilde \lin v:= \Delta v + |\tilde A|^2 v - \tilde X \cdot \nabla v + v $ associated to normal perturbations of $\tilde H + \tilde X \cdot \nu$.  

\begin{proposition} 
\label{prop:lin-cap}
Given $ E \in C^{0, \alpha}(\tilde \ccal)$, there exist a function $ v \in C^{2,\alpha}(\tilde \ccal)$ and a constant $C$ such that 
	\begin{gather*}
	\tilde \lin_{\tilde \ccal} v = E, \qquad \left. v \right|_{\pd \tilde \ccal} = 0,\\
	\|  v \|_{C^{2, \alpha}} \leq C \|  E\|_{C^{0,\alpha}}.
	\end{gather*}
\end{proposition} 

\begin{proof}
Let $\mathbf S^2$ be the standard 2-sphere and $S$ be the sphere of radius $\sqrt 2$, both equipped with the metrics induced by their respective embeddings into $E^3$.  The linear operator of interest on $S$ is $\tilde \lin_S = \Delta_{S} + |A_{S}|^2 + 1 = \frac{1}{2}  \Delta_{\mathbf S^2} +2 $.  The existence of a unique solution satisfying the estimate above is standard on a hemisphere of $S$ thanks to the study of eigenvalues of the Laplace operator on the unit sphere $\mathbf S^2$ (see for example \cite{stein-weiss;introduction-fourier}).  We obtain the result for $\tilde \lin_{\tilde \ccal}$ by treating $\tilde \ccal$ as a perturbation of a hemisphere of $S$.
\end{proof}

\subsection{The linearized equation on the inner disk $\tilde \dcal$}

The existence of a solution for the Dirichlet problem $\tilde \lin  v =  E$, $\left. v \right|_{\pd \tilde \dcal} =0$ with estimates similar to the ones in Proposition \ref{prop:lin-cap} follows from standard theory in PDEs. 


\subsection{The linearized equation on the outer plane $\tilde \pcal$}
\label{ssec:linear-P}


Let $\bar R = \tilde R e^{\tau(a+\ubar a)/\tilde R }$. We denote by $B_{\bar R} \subset \R^2$ the disk of radius $\bar R$ centered at the origin and by $\Omega:= \R^2\setminus B_{\bar R}$ the plane with the disk of radius $\bar R$ removed. Because $\tilde \pcal$ only differs from $\Omega$ in a small neighborhood of the boundary, any results and estimates we obtain for the solution to the linearized equation $\tilde \lin_{\Omega} v =E$ are also valid for the solution to the linearized equation on $\tilde \pcal$ by using perturbation theory.  

The preliminary estimates from \cite{mine;part2} show that if the nonhomogeneous term $E$ decays as $r^{-1}$, where $r$ is the distance to the origin, then the solution $v$  of the linearized equation on $\Omega$ has linear growth, bounded gradient, and a Hessian decaying as $r^{-1}$. We emphasize the fact that $v$ is not bounded but is asymptotic to a cone at infinity. The weighted norms below and in Definition \ref{def:C-cone} capture this information. 

Let $\xi$ be a point in a region $\ncal$.  For $v \in C^{r,\alpha}_{loc}(\ncal)$, $r=0,2$, we define the following norms: 
	\[
	\|  v :C_{\ast}^{r,\alpha}(\ncal)\|=\max\left(\max_{0\leq j\leq r} \| D^jv(\xi) |\xi|^{-r+1+j} \|_{C^0(\ncal)},\sup_{\xi\in \ncal} \left([D^rv]_{\alpha, B(\xi)\cap \ncal} |\xi|^{1+\alpha}\right) \right),
	\]
where $B(\xi)$ denotes the geodesic ball of radius $1$ centered at $\xi$, and $[v]_{\alpha, B}$ is the usual H\"older semi-norm 
	\[
	[v]_{\alpha, B} = \sup_{\eta, \eta' \in B} \frac{ |v(\eta) - v(\eta')|}{|\eta-\eta'|^{\alpha}}.
	\]
\begin{definition}
\label{def:C-ast}
$C_{\ast}^{r,\alpha}(\Omega)$ is the space of functions in $C^{r,\alpha}_{loc}(\Omega)$ with finite  $C_{\ast}^{r,\alpha}$ norm and whose graphs over $\Omega$ satisfy the imposed symmetries (ii) and (iii)  from Proposition \ref{prop:smooth-initial-surfaces}. 
\end{definition}
Note that the dependence of $C_{\ast}^{r,\alpha}(\Omega)$ on $m$ is implicit here and in the rest of the article.

\begin{proposition}
\label{prop:lin-plane} Given $E\in C_{\ast}^{0,\alpha}(\Omega)$, there exist a unique $v\in C_{\ast}^{2,\alpha}(\Omega)$ and a constant $C$ depending only on $\bar R$ so that 
	\begin{gather}
	\notag \tilde \lin v = \Delta v - \xi\cdot \nabla v + v  = E, \qquad v|_{\pd \Omega}=0,\\
	\label{eq:linear-P}\| v: C_{\ast}^{2,\alpha}(\Omega) \|  \leq C \|E:C_{\ast}^{0,\alpha}(\Omega)\|.
	\end{gather}
\end{proposition}

\begin{proof}
In Lemma 5 and Theorem 7 of \cite{mine;part2}, we showed the existence and uniqueness of a weak solution $v$ (note that the roles of $u$'s and $v$'s are swapped in the mentioned article). The operator is elliptic, so a weak solution is also a strong smooth solution in $\Omega$. 

We now prove the estimate on $\| v:C^{2,\alpha}_{\ast}(\Omega)\|$. Let us first recall how to obtain the bounds on $|v|$. Since  $1-2\bar R^2 <0$, the function $v_k=k(r-\frac{\bar R^2}{r})$ with $r=|\xi|$ is a supersolution and satisfies
	\[ 
	\tilde \lin v_k=-k\frac{\bar R^2}{r^3}+\frac{k}{r}(1-2 \bar R^2) \leq E.
	\]
for  $k \geq \sup_{\xi \in \Omega} (|E(\xi)| |\xi| /(1-2\bar R^2))$. On the sector $\Omega_{\bar R, m}:=\{ (r \cos \theta, r\sin \theta) \in \Omega \mid \theta \in(-\pi/m, \pi/m)\}$, we have $\tilde \lin (v_k - v) \leq 0$ and $v_k-v \geq 0$ on $\pd \Omega_{\bar R,m}$. Using the symmetries and the maximum principle on a sector (Theorem 7 \cite{mine;part2}), we get 	
 	\begin{equation}
	\label{eq:linear-growth}
	|v(\xi)|\leq k |\xi|.
	\end{equation}

Using the standard theory of functions in H\"older spaces, we can extend a function $v\in C^{2,\alpha}(\Omega)$ to a function $v'$ in $C^{2,\alpha}(\R^2)$ in such a way that 
	\[
	\| v':C^{2,\alpha}(B_{\bar R+1})\| \leq C \| v:C^{2,\alpha}(B_{\bar R+1} \cap \Omega)\|,
	\]
where the constant $C$ is independent of $v$. Let us fix once and for all such an extension. We define $E'$ to be $\tilde \lin v'$. The advantage of working with functions on $\R^2$ is that the heat equation \eqref{eq:heat} below is now defined on the whole space. We can then use the well-known heat kernel for $\R^2$. 

Before tackling the heat equation, we need to show that 
	\begin{equation}
	\label{eq:E'-norm}
	\|E':C^{0,\alpha}_{\ast}(\R^2)\| \leq C \| E: C^{0,\alpha}_{\ast}(\Omega)\|.
	\end{equation}
By definition, we have $E'=E$ on $\Omega$ and 
	\[
	\|E':C^{0,\alpha}(B_{\bar R+1})\| \leq C \|v':C^{2,\alpha}(B_{\bar R+1})\| \leq C \|v:C^{2,\alpha}(B_{\bar R+1} \cap \Omega)\|).
	\]
Because $v=0$ on $\pd \Omega$, boundary estimates for linear elliptic equations (see Lemmas 6.4 and 6.5 in \cite{gilbarg-trudinger;elliptic-equations} for example) give a $\delta>0$ for which 
	\begin{equation}
	\label{eq:v-boundary}
	\| v :C^{2,\alpha}(B_{\bar R + \delta} \cap \Omega)\| \leq C (\| v : C^0(B_{\bar R+1} \cap \Omega)\| + \| E :C^{0,\alpha}(B_{\bar R+1} \cap \Omega) \|). 	\end{equation}
By Schauder interior estimates, we have 
	\begin{equation}
	\label{eq:v-interior}
	\| v : C^{2,\alpha} (B_{\bar R+1}\setminus B_{\bar R + \delta})\| \leq C (\| v : C^0(B_{\bar R+1+\delta} \cap \Omega)\| + \| E :C^{0,\alpha}(B_{\bar R+1+\delta} \cap \Omega) \|,
	\end{equation}
where $C$ depends on $\delta$. 
Combining the inequalities \eqref{eq:linear-growth}, \eqref{eq:v-boundary}, and \eqref{eq:v-interior}, we obtain \eqref{eq:E'-norm}.

Using the change of variables $\xi=\eta/\sqrt{2(1-t)}$, we define the new function  
	$
	u(t,\eta) := \sqrt{2(1-t)} \ v'\left(\frac{\eta}{\sqrt{2(1-t)}}\right)
	$
which satisfies the heat equation 
	\begin{align}
	\label{eq:heat}
	\pd_t u -\Delta u 
	= - \frac{1}{\sqrt{2(1-t)}} E'\left(\frac{\eta}{ \sqrt{2(1-t)}}\right) 
	\end{align}
on the parabolic cylinder $(0,1) \times \R^2 $.
Because of the scaling, to prove \eqref{eq:linear-P} it suffices to show that for $|\eta|=\bar R$, we have
	\[
	 \sup_{t \in (\frac{1}{2}, 1)} |\pd_{\eta_i \eta_j} u(t, \eta)| + \sup_{t \in (\frac{1}{2}, 1)}[\pd_{\eta_i \eta_j} u(t, \eta)]_{\alpha,B(\eta)} \leq \|E: C^{0,\alpha}_{\ast}(\Omega)\|.
	\]
The fundamental solution of the heat equation $\pd_t u -\Delta u =f$, $u|_{t=0}=u_0$ is given by 
	\[
	G(t,\eta) = ({4\pi t})^{-1} e^{-\frac{|\eta|^2}{4t}}
	\]
and 
	\begin{multline}
	\label{eq:pd2etaG}
	\pd^2_{ \eta_i \pd \eta_j} u(t,\eta) = \int_{0}^t d\tau \int_{\R^2}\pd^2_{ \eta_i \pd \eta_j}G(t-\tau, \eta-\zeta) [f(\tau,\zeta) - f(\tau, \eta)] d\zeta\\
	+  \int_{\R^2} \pd^2_{ \eta_i \pd \eta_j}G(t, \eta-\zeta) u_0(\zeta) d\zeta.
	\end{multline}
(see \cite{eidelman;parabolic-systems} pp17-20 for example.)

Many of the regularity theorems in the literature involve both the temporal and spatial derivatives of $u$. Here, we only need an estimate on the spatial H\"older semi-norm of $\pd^2 u$ for fixed time $t$. For this, we treat the terms on the right hand side of \eqref{eq:pd2etaG} separately.  Lemma \ref{lem:heat-seminorm} below is a classical result about solutions to the heat equation with zero initial condition and can be found in \cite{ladyzhenskaja;linear-parabolic-type} p.275. We provide an alternate proof here.

\begin{lemma}[A characterization of $C^{\alpha}$]
\label{lem:Calpha}
A function $w: \R^n \to \R$ is $C^{\alpha}$ if and only if for any $\eps>0$, there is an $w_{\eps} \in C^1$ such that 	
	\[
	\| w - w_{\eps} \|_{L^{\infty}} \leq C_0 \eps^{\alpha}, \quad \| \nabla w_{\eps} \|_{L^{\infty}} \leq C_1 \eps^{\alpha-1},
	\]
with $C_0$ and $C_1$ dependent on $w$ but not on $\eps$. 
\end{lemma}

\begin{proof}
Let $w \in C^{\alpha}$. We  take $w_{\eps}$ to be the convolution $ w \ast \varphi_{\eps}$ where $\{ \varphi_{\eps}\}$ is a family of smooth functions such that $\varphi_1 (x)$ is compactly supported, $0\leq \varphi_1(x) \leq 1$, $\int \varphi_1 =1$, and $\varphi_{\eps} (x) = \eps ^{-1} \varphi_1(x/\eps)$.

Conversely, given $x$ and $y$, we can pick $\eps = |x-y|$ and obtain
	\[
	|w(x) -w(y)| \leq |w(x) - w_{\eps}(x)| + |w_{\eps}(x) - w_{\eps}(y)| + |w_{\eps}(y) -w(y)| \leq (2C_0+C_1) \eps^{\alpha},
	\]
with the estimates 
	\[
	\| w\|_{L^{\infty}} \leq C_0 \eps^{\alpha} + \| w_{\eps}\|_{L^{\infty}}, \quad [w]_{\alpha} \leq 2C_0+C_1. \qedhere
	\]
\end{proof}

\begin{lemma}
\label{lem:heat-seminorm}
Let $f$ be a function in $C_{loc}^{0, \alpha}((0, 1) \times \R^2)$ and define
	\[
	\bar u (t,\eta)=  \int_{0}^t d\tau \int_{\R^2}G(t-\tau, \eta-\zeta) f(\tau,\zeta) d\zeta.
	\]
We have
	\begin{gather*}
	\sup_{0<t<T} \left[\pd^2_{ \eta_i \pd \eta_j} \bar u(t,\cdot)\right]_{\alpha} \leq C \sup_{0<t<T} [f(t,\cdot)]_{\alpha},
	\\
	\sup_{0<t<T} \left\|\pd^2_{ \eta_i \pd \eta_j} \bar u(t,\cdot)\right\|_{L^{\infty}} \leq C \sup_{0<t<T} [f(t,\cdot)]_{\alpha}+ C \sup_{0<t<T} \| f(t,\cdot) \|_{L^{\infty}}.
	\end{gather*}
\end{lemma}

\begin{proof} 
Fix a time $t$ in $(0,T)$ and let $\eps$ be an arbitrary constant satisfying $0<\eps<\sqrt t$. In this proof, we will use the notation $[f]_{T,\alpha} := \sup_{0<t<T} [f(t,\cdot)]_{\alpha}$. We have
	\begin{multline*}
	\pd^2_{\eta_i  \eta_j} \bar u(t,\eta) =
	 \int_{t-\eps^2}^t d\tau \int_{\R^2} \pd^2_{\eta_i \eta_j}G(t-\tau, \eta-\zeta) [f(\tau,\zeta)- f(\tau, \eta)] d\zeta \\
	+ \int_0^{t-\eps^2} d\tau \int_{\R^2} \pd^2_{\eta_i \eta_j}G(t-\tau, \eta-\zeta) [f(\tau,\zeta) - f(\tau, \eta)] d\zeta.
	\end{multline*}
Define  $g_{\eps}$ to be the first term and $w_{\eps}$ to be second term on the right hand side. Note that the function $w_{\eps}(t,\cdot):\R^2 \to \R$ is continuously differentiable. 

We will prove that $w:=\pd^2_{\eta_i  \eta_j} \bar u=g_{\eps}+w_{\eps}$ is H\"older continuous in the variable $\eta$ using Lemma \ref{lem:Calpha}. After performing the change of variables $s=t-\tau$ and $y=\eta-\zeta$, we obtain
	\begin{align*}
	|g_\eps| &\leq \int_{0}^{\eps^2} \int \left| G(s,y) \left( \frac{y_iy_j}{4s^2} - \frac{\delta_{ij}}{2s} \right) [f]_{T,\alpha}|y|^{\alpha}\right|dy ds\\
			& = \int_0^{\eps^2}\int \left|\frac{1}{({4\pi s})}  e^{\frac{-|y|^2}{4s}} \left( \frac{y_iy_j}{4s^2} - \frac{\delta_{ij}}{2s} \right) [f]_{T,\alpha}|y|^{\alpha}\right|dy ds.
	\end{align*}
With $z = \frac{y}{2\sqrt s}$, the inequality above becomes
	\[
	|g_{\eps}| \leq C[f]_{T,\alpha}  \int_0^{\eps^2}\int  \frac{|z|^2+1}{s^{1-\alpha/2}} |z|^{\alpha} e^{-|z|^2} dz ds \leq C[f]_{T,\alpha} \eps^{\alpha}.
	\]
We now estimate $|\pd_{\eta_k} w_{\eps}|$. Since the variable $\tau$ stays away from $t$, the integrals in $w_{\eps}$ converge. Note also that $\int \pd^2_{\eta_i \eta_j}G(t-\tau, \eta-\zeta) d\zeta=0$ so we can write 
	\[
	\pd_{\eta_k}w_{\eps} = \int_0^{t-\eps^2} d\tau \int \pd^3_{\eta_i \eta_j \eta_k}G(t-\tau, \eta-\zeta) [f(\tau,\zeta) - f(\tau, \eta)]  d\zeta.
	\]
Performing the same changes of variables above, we get  
	\begin{align*}
	|\pd_{\eta_k}w_{\eps}| &\leq [f]_{T,\alpha}\int_{\eps^2}^t \int  \left|G(s,y) \left( -\frac{y_iy_jy_k}{8s^3} + 				 			  \frac{\delta_{ij}y_k+\delta_{jk}y_i+\delta_{ki}y_j}{4s^2} \right) |y|^{\alpha} \right|dy ds \\
				 &\leq C [f]_{T,\alpha} \int_{\eps^2}^t \int \frac{|z|^3+1}{s^{(3-\alpha)/2}} |z|^{\alpha} e^{-|z|^2} dz ds 
				\leq C[f]_{T,\alpha} \eps ^{\alpha-1},
	\end{align*}
where we used $0<\eps^2<t$ in the last inequality. We also have
	\begin{align*}
	| w_{\eps} |&\leq \int_0^{t-\eps^2} d\tau \int \left| \pd^2_{x_i x_j}G(t-\tau, x-\xi) [f(\tau,\xi) -f(\tau, x)] \right| d\xi\\
			& \leq  C[f]_{T,\alpha} \int_{\eps^2}^t\int  \frac{|\eta|^2+1}{s^{1-\alpha/2}} |\eta|^{\alpha} e^{-|\eta|^2} d\eta ds  \\
			&\leq  C[f]_{T,\alpha} t^{\alpha/2}. 
	\end{align*}
If $|\eta-\eta'|<\sqrt t$ for $\eta, \eta' \in \R^2$, we can choose $\eps=|\eta-\eta'|$ and the proof of Lemma \ref{lem:Calpha} gives $|w(t,\eta) - w(t,\eta')| \leq C [f]_{T,\alpha} |\eta-\eta'|^{\alpha}$. If $|\eta-\eta'| \geq \sqrt t$, we choose $\eps=\sqrt t/2$ and get 
	\begin{align*}
	 |w(t,\eta) & - w(t,\eta')| \\&\leq |w(t,\eta) - w_{\eps}(t,\eta)| + |w_{\eps}(t, \eta) - w_{\eps}(t,\eta')| + |w_{\eps}(t,\eta')-w(t,\eta')| \\
	& \leq C  [f]_{T,\alpha} \eps^{\alpha} + 2 \| w_{\eps} (t, \cdot)\|_{L^{\infty}(\R^2)} |\eta-\eta'|^{\alpha} t^{-\alpha/2}\\
	& \leq C [f]_{T,\alpha} |\eta-\eta'|^{\alpha}.
	\end{align*}
It is clear that $\|w(t,\cdot)\|_{L^{\infty}(\R^2)} \leq C[f]_{T,\alpha}$. 
\end{proof}
It remains to estimate the last term in \eqref{eq:pd2etaG} for $t \in (\frac{1}{2},1)$, $|\eta|=\bar R$, with $u_0=\sqrt 2 v'(\eta/\sqrt 2)$. Let us denote it by $h$. We have
	\begin{align*}
	|h(t,\eta)|&\leq  \int_{\R^2} \left| G(t,y) \left( \frac{y_iy_j}{4t^2} - \frac{\delta_{ij}}{2t} \right) u_0(\eta-y) \right| dy \\
			&\leq C  \|E:C^{0,\alpha}_{\ast}(\Omega)\| \int_{\R^2} t^{-2}  e^{-C' \frac{|y|^2}{t}} |\eta-y| dy\\
			&\leq C  \|E:C^{0,\alpha}_{\ast}(\Omega)\|,
	\end{align*}
where we used the linear growth of $v$ given by \eqref{eq:linear-growth} for the second line. For the H\"older semi-norm, the Mean Value Theorem implies
	\begin{align*}
	|h(t,\eta) &- h(t,\eta')| \leq \int_{\R^2} \int_0^1 \left| \pd^3_{\eta_k \eta_i \eta_j} G(t, \eta'-\zeta+s(\eta-\eta'))\right| |\eta_k-\eta'_k| |u_0(\zeta)| ds d\zeta \\	
		&\leq C  \|E:C^{0,\alpha}_{\ast}(\Omega)\| |\eta-\eta'|^{\alpha} \int_{\R^2} \int_0^1 t^{-5/2} e^{-C'' \frac{|\eta'-\zeta+s(\eta-\eta')|^2}{t}} |\zeta| ds d\zeta\\
		&\leq C  \|E:C^{0,\alpha}_{\ast}(\Omega)\| |\eta-\eta'|^{\alpha} \left( \int_{B_{10}(0)}  |\zeta| d\zeta + \int_{\R^2 \setminus B_{10}(0)} e^{- C''\frac{|\zeta|^2}{t}} |\zeta| d\zeta \right)\\
		& \leq C  \|E:C^{0,\alpha}_{\ast}(\Omega)\| |\eta-\eta'|^{\alpha}, 
	\end{align*}
where  the values of the constants $C$ and $C''$ may be adjusted in each line. 

We conclude the proof of Proposition \ref{prop:lin-plane} by applying Lemmas \ref{lem:Calpha} and \ref{lem:heat-seminorm} and noting that for $t\in(0,1)$, $ [f(t,\cdot)]_{\alpha}$ and $\| f(t,\cdot) \|_{L^{\infty}}$ are both bounded by $C  \|E:C^{0,\alpha}_{\ast}(\Omega)\|$.	
\end{proof}

Unfortunately, $C^{2,\alpha}_{\ast}(\Omega)$ is not suitable for the Schauder Fixed Point Theorem in Section \ref{sec:fpt} because bounded sets in $C^{2,\alpha}_{\ast}(\Omega)$ are not compact in $C^{2,\alpha'}_{\ast}(\Omega)$, $0<\alpha'<\alpha<1$. The self-shrinkers we construct are not asymptotically planar but tend to cones at infinity. We take advantage of this asymptotic behavior in the definition below.

\begin{definition} 
\label{def:C-cone}
We define $C^{2,\alpha}_{cone}(\Omega)$ to be the space of functions $v$ in $C^{2,\alpha}_{loc}(\Omega)$ that satisfy the following requirements:
\begin{enumerate}
\item there are functions $\varphi: S^1 \to \R$ and $w: \Omega \to \R$ such that 
	\[
	v (\xi) = \varphi(\xi/|\xi|) |\xi| + w(\xi),
	\]
\item $ \| \varphi: C^{2,\alpha} (S^1) \| < \infty$,
\item $ \| w : C^{2,\alpha}(\Omega, |\xi|^{-1})\|< \infty$,
\item  $\| D^2 w: C^{0,\alpha}_{\ast} (\Omega) \| < \infty$ and $\| \xi \cdot \nabla w : C^{0,\alpha}_{\ast}(\Omega) \| < \infty$.
\end{enumerate}
The $C^{2,\alpha}_{cone}(\Omega)$ norm of $v$ is the maximum of the quantities in (ii), (iii), and (iv). 
\end{definition} 
It is easy to see that if a decomposition of $v$ into $\varphi$ and $w$ exists, it is unique. The $C^{2,\alpha}_{cone}(\Omega)$-norm is therefore well defined. Thanks to the compactness of $S^1$ and the decay rate in (iii)(iv), bounded sets in $C^{2,\alpha}_{cone}(\Omega)$ are compact in $C^{2,\alpha'}_{cone}(\Omega)$ for $0<\alpha'<\alpha<1$.  

The requirements in (iv) of Definition \ref{def:C-cone} were added to ensure that $\tilde \lin(C^{2,\alpha}_{cone}) \subseteq C^{0,\alpha}_{\ast}(\Omega)$. The fact that $w \in C^{0,\alpha}_{\ast}(\Omega)$ comes from the uniform bound on $||\xi| w| $ and  the estimate $[w]_{2+\alpha, B(\xi) \cap \Omega} \leq C [D^2 w]_{\alpha, B(\xi )\cap \Omega}$, where $C$ depends on $\alpha$ only. 
\begin{lemma}
\label{lem:cone-op}
Given $E \in C^{0,\alpha}_{\ast}(\Omega)$, the solution $v$ to $\tilde \lin v =E$, $v|_{\pd \Omega} =0$ given in Proposition \ref{prop:lin-plane} satisfies 
	\[
	\| v : C^{2,\alpha}_{cone}(\Omega) \| \leq C \| E : C^{0,\alpha}_{\ast} (\Omega) \|.
	\]
\end{lemma}
	
\begin{proof}	
Let us denote by $f$ the function $E- \Delta v$ and by $K$ the constant $\| E : C^{0,\alpha}_{\ast} (\Omega) \|$. We have  $\| f : C^{0,\alpha}_{\ast}(\Omega)\| \leq C K$ by Proposition \ref{prop:lin-plane}. We will use polar coordinates $(r, \theta)$ in $\Omega$ and abuse notations by sometimes writing $v(r,\theta)$ to mean $v(r \cos \theta, r \sin \theta) = v(\xi)$. 

Consider the linear first order equation $- \xi \cdot \nabla v + v = f$, which can be rewritten in polar coordinates as
	$
	- r \pd_r v + v = f.
	$
For fixed $\theta$, the solution is
	\[
	v (r,\theta) = r  \int_r^{\infty} \frac{f(s,\theta)}{s^2}\ ds + c_1 r, 
	\]
where $c_1$ is a function of $\theta$ only. From the boundary condition $v|_{\pd \Omega}=0$, we obtain $c_1 = - \int_R^{\infty}\frac{f(s,\theta)}{s^2}\ ds$. Define 
	\begin{equation*}
	w (r, \theta) := r  \int_r^{\infty} \frac{f(s,\theta)}{s^2}\ ds \quad \textrm{and} \quad \varphi(\theta) : = - \int_R^{\infty}\frac{f(s,\theta)}{s^2}\ ds.
	\end{equation*}
We know the integrals above exist because $|f(s,\theta)| < \frac{CK}{s}$. Moreover, 
	\begin{equation}
	\label{eq:w-C0}
	|w(r,\theta)| \leq r \int_{r}^{\infty} \frac{|f|}{s^2} ds \leq \frac{CK}{r}
	\end{equation} 
so $v (r,\theta) \to \varphi(\theta) r$ uniformly as $r \to \infty$. 

For $\lambda \in [1, \infty)$ and $ \xi \in \Omega$, we define the scaled functions 
	\[
	v_{\lambda} (\xi) = \lambda^{-1} v (\lambda \xi),
	\]
which satisfy $\| v_{\lambda} : C^{2,\alpha}_{cone}(\Omega) \| = \| v : C^{2,\alpha}_{cone}(\Omega)\|$. For a fixed ball $B_j$, the $v _{\lambda}$'s are bounded in $C^{2,\alpha}(B_j \cap \Omega)$ and by the Arzela-Ascoli Theorem, given $\alpha'<\alpha$, there exist a function $v_{\infty} \in C^{2,\alpha}(B_j \cap \Omega)$ and a subsequence $v_{\lambda_k}$ that converges  in $C^{2, \alpha'}(B_j \cap \Omega)$ to $v_{\infty}$. From the uniqueness of the limit, we have $v_{\infty} = \varphi(\xi/|\xi|) |\xi|$. The fact that $ v_{\infty} \in C^{2,\alpha}(B_j\cap \Omega)$ and the bound on  $v_{\lambda_k}$ imply 
	$
	\| \varphi : C^{2,\alpha}(S^1) \| \leq C K.
	$
With straightforward computations, one can show that 
	\[
	\| D^2[\varphi(\xi/|\xi|) |\xi|] : C^{0,\alpha}_{\ast} (\Omega) \| \leq C K;
	\]
therefore 
	$
	\| D^2 w : C^{0,\alpha}_{\ast}(\Omega) \| =\| D^2 [v - \varphi(\xi/|\xi|) |\xi|] : C^{0,\alpha}_{\ast}(\Omega) \| \leq CK.
	$
This last estimate and \eqref{eq:w-C0} give us that $w \in C^{0,\alpha}_{\ast}(\Omega)$. Recall that
	\[
	\tilde L v = \Delta (r \varphi) + \Delta w - r \pd_r w + w = E.
	\]
Hence, $r \pd_r w = \Delta ( r \varphi) + \Delta w  + w -E \in C^{0,\alpha}_{\ast}(\Omega)$ and its norm is bounded by $CK$. 

We now finish the proof by showing that $|Dw| \leq C K |\xi|^{-1}$. Without loss of generality, we can assume that $|\xi| > 10$. In particular, this means that the ball of radius $2$ centered at $\xi$, $B_2(\xi)$, is in $\Omega$. The function $w$ satisfies the Poisson equation $\Delta w = F$ in $B_2(\xi)$, where $F = E -w + r \pd_r w - \Delta (r \varphi)$. Interior H\"older estimates for solutions to Poisson's equation give us the desired bound on $|Dw|$ (see Lemma 4.6 in \cite{gilbarg-trudinger;elliptic-equations} with $R=1$). We outline the relevant part of the proof from \cite{gilbarg-trudinger;elliptic-equations} in the next paragraph.

We can write $w = w'+ w''$ where $w'$ is a harmonic function on $B_2(\xi)$ (with boundary conditions $w'=w$ on $\pd B_2(\xi)$) and $w''$ is the Newtonian potential of $F$ in $B_2(\xi)$.  Let $\Gamma(\xi-\xi')  = \frac{1}{2\pi} \log|\xi-\xi'|$ be the fundamental solution of the Laplace equation. Standard theory on the Laplace operator  gives
	\[
	\|Dw'\|_{C^{0}(B_1(\xi))}  \leq C \sup_{B_2(\xi)} |w'| \leq C \sup_{\pd B_2(\xi)} |w'| \leq C K |\xi|^{-1}
	\]
and 
	\[
	D_i w'' (x) = \int_{B_2(\xi)} D_i \Gamma (x-y) F(y) dy, i=1,2.
	\]
Combining the above equation with $| D \Gamma (\xi - \xi') | \leq C |\xi-\xi'|^{-1}$ and $| F | < CK |\xi|^{-1}$, we obtain the estimate on $| Dw|$. 
\end{proof}

\subsection{The linearized equation on $\tilde M$}

Once the correct Banach spaces of functions are defined, the rest of the construction (solving the linearized equation $\tilde \lin_{\tilde M} v =E$ on $\tilde M$ and using a Fixed Point Theorem for the solution to the nonlinear equation \eqref{eq:self-shrinker}) follows the same lines as in \cite{kapouleas;embedded-minimal-surfaces} or \cite{mine;self-trans}. We provide the few finishing touches here to give a coherent ending to this article.

We define a global norm on $\tilde M$ from the various norms used on $\Sigma$, $\tilde \dcal$, $\tilde \pcal$ by essentially taking the maximum of all these norms. The factor $e^{-5 \delta_s /\tau}$ takes into account that our functions are decaying on the overlapping regions and the factor $\tau^{10}$ reflects a loss in exponential decay incurred because we scale and cut the local solutions in the proof of Theorem \ref{thm:linear-invertible} below. It is not significant compared to the exponential weight. 

Let us recall that $\hcal$ is the homothety of ratio $\tau$ centered at the origin.
\begin{definition}
Given $ v \in C^{0,\alpha}_{loc}(\tilde M)$, we define $\| v \|_0$ to be the maximum of the quantities below, where $b_0 = e^{-5 \delta_s /\tau}$,

	\begin{enumerate}
	\item $
		\tau \|  v \circ \hcal: C^{0,\alpha}(M \cap (\Sigma \cup \ccal \cup \dcal), g_{M}, \max(e^{-\gamma s}, b_0)) \|,
		$
	\item 
		$
		b_0^{-1} \|  v: C^{0,\alpha}_{\ast}(\tilde \pcal\setminus \tilde \Sigma) \|.
		$
	\end{enumerate}
	Given $ v \in C^{2,\alpha}_{loc}(\tilde M)$, we define $\| v \|_2$ to be the maximum of the quantities below, where $b_2= b_0/\tau^{10}$,

	\begin{enumerate}
	\item $
		\tau^{-1} \|  v \circ \hcal: C^{2,\alpha}(M \cap (\Sigma \cup \ccal \cup \dcal), g_{M}, \max(e^{-\gamma s}, b_2)) \|,
		$
	\item 
		$
		b_2^{-1} \|  v: C^{2,\alpha}_{cone}(\tilde \pcal\setminus \tilde \Sigma) \|.
		$
	\end{enumerate}
\end{definition} 
Note that $\hcal^{\ast} g_{\tilde M} = \tau^2 g_M$. For any  function $ v\in C^{2,\alpha}_{loc}(\tilde M)$ supported on $\tilde \Sigma \cup \tilde \ccal \cup \tilde \dcal$,  the corresponding function $\bar v:= \tau^{-1}  v \circ \hcal$ has the following property
	\[
	\|  v \|_2 = \| \bar v : C^{2,\alpha}(M \cap (\Sigma \cup \ccal \cup \dcal), g_{M}, e^{-\gamma s}) \|. 
	\]
Similarly, taking $\bar E:= \tau E \circ \hcal$ for a function $E$ supported on $\tilde \Sigma \cup \tilde \ccal \cup \tilde \dcal$ gives
	\[
	\| E\|_0 = \| \bar E: C^{0,\alpha}(M \cap (\Sigma \cup \ccal \cup \dcal), g_{M}, e^{-\gamma s}) \|. 
	\]
Moreover, these new functions $\bar v$ and $\bar E$ satisfy $\lin_M \bar v = \bar E$ if and only if $\tilde \lin_{\tilde M} v =  E$. 

To simplify the notations later on, we define the linear map $\Theta: [-\zeta \tau, \zeta \tau] \to C^{\infty}_{loc} (\tilde M)$ by 
	\[
	\Theta(b) = \tau^{-1} \hcal_{\ast} (b w).
	\]

\begin{theorem}
\label{thm:linear-invertible}
Given $ E \in C^{0,\alpha}_{loc}(\tilde M)$ with finite norm $\|  E \|_0$, there exist $ v_E \in C^{2,\alpha}_{loc}(\tilde M)$ and a constant $b_E$  uniquely determined by the construction below, such that
	\[
	\tilde \lin_{\tilde M}  v_E =  E + \Theta(b_E),
	\]
and 
	\[
	\|  v_E\|_2 \leq C \|  E\|_0, \quad |b_E |\leq C \|  E\|_0.
	\]
\end{theorem}

\begin{proof}
Let  $\psi$ be the cut-off function on $\tilde M$ defined by 
	$
	\psi := \psi[5 \delta_s /\tau, 5 \delta_s /\tau -1]\circ s$ on $\tilde \Sigma$ and $\psi\equiv 0$ on the rest of $\tilde M$. 

We take $E_0:= E$ and proceed by induction; given $ E_{n-1}$, we define $ E_n$, $ v_n$ and $b_n$ in the following way.  First, we apply Proposition \ref{prop:7-1} on the desingularizing piece $\Sigma = \Sigma[\tilde R (b) , b, \tau]$ with $E' = \tau (\psi  E_{n-1}) \circ \hcal$  to  obtain  $v_{E'}$ and $b_{E'}$. We take $b_n:=b_{E'}$ and define the  function $ u' := \tau \hcal_{\ast} v_{E'}$ which satisfies
	\begin{gather*}
	\tilde \lin_{\tilde M}  u' = \psi  E_{n-1} + \Theta(b_n),	\end{gather*}
The function $\psi u'$  can be extended smoothly by zero to the rest of $\tilde M$ and from the estimate in Proposition \ref{prop:7-1}, we have
	\begin{equation}
	\label{eq:est-u'}
	\| \psi  u' \|_2 \leq C \| E_{n-1}\|_0.
	\end{equation}
Note for the next steps that $
	\tilde \lin_{\tilde M} (\psi u') =  \psi^2  E_{n-1} +   [\tilde \lin, \psi]  u' + \Theta(b_n), 
	$
where we used the notation $[\tilde \lin_{\tilde M}, \psi]  u':= \tilde \lin_{\tilde M} (\psi) u' - \psi(\tilde \lin_{\tilde M}  u')$. 

The function $ E'':=  E_{n-1}  - \psi^2  E_{n-1} -   [\tilde \lin, \psi]  u'$ is supported on $s \geq \frac{5 \delta_s}{\tau} -1$, therefore it can be decomposed into $ E'' =  E''_{\ccal} +  E''_{\dcal}+ E''_{\pcal}$ where each $ E''_{\ncal}$ is supported in $\tilde \ncal$. From the discussion in Section \ref{sec:linear}, there exist functions $ u''_{\ccal}$, $ u''_{\dcal}$, and $ u''_{\pcal}$ that satisfy for $\ncal =\ccal, \dcal, \pcal$,
	\begin{gather*}
	\tilde \lin  u''_{\ncal} =  E''_{\ncal},  \textrm{ in } \tilde  \ncal,\\
	 u''_{\ncal} = 0 \textrm{ on } \pd \tilde \ncal.
	\end{gather*} 
Let $ u''$ be the continous function $ u''_{\ccal} + u''_{\dcal}+ u''_{\pcal}$ extended by zero to the rest of $\tilde M$. On each of the bounded pieces $\tilde \ccal $, $\tilde \dcal$, and $\tilde \pcal \cap \tilde \Sigma$ we have 
	\begin{align}
	\notag \|  u'': C^{2,\alpha} \| &\leq C \|  E_{n-1}  - \psi^2  E_{n-1} -   [\tilde \lin, \psi]  u' : C^{0,\alpha}\|,\\
	\label{eq:est-u''}					&\leq C \tau^{-1-\alpha}e^{- \gamma  (5\delta_s/\tau -1)} \|  E_{n-1} \|_0,
		\end{align}
and on $\tilde \pcal$,
	\begin{equation}
	\label{eq:u''pcal}
	\|  u'': C^{2,\alpha}_{cone}(\tilde \pcal) \| \leq C \tau^{-1-\alpha}e^{- \gamma  (5\delta_s/\tau -1)} \|  E_{n-1} \|_0.
	\end{equation}
	
We define another cut-off function $\psi':= \psi[\ubar a, \ubar a+1]\circ s $ on $\tilde M$ and a function  $ v_n = \psi u' +\psi' u''$.  Let us recall that $\ubar a = 8| \log \tau|$, as in Definition \ref{def:components}. Its logarithmic dependence on $\tau$ will be crucial for proving that we have a contraction \eqref{eq:contraction}. 
	Since  $\psi'\equiv 1$ on the supports of $1-\psi^2$ and $[\tilde \lin, \psi]$, $ v_n$ satisfies 
	\begin{gather*}
	\tilde \lin  v_n =   E_{n-1}  +   [\tilde \lin, \psi']  u''+ \Theta(b_n), \\
	\|  v_n \|_2 \leq C \|  E_{n-1} \|_0,
	\end{gather*}
where the inequality follows from \eqref{eq:est-u'}, \eqref{eq:est-u''}, and \eqref{eq:u''pcal}. 
We define   $ E_n = - [\tilde \lin, \psi'] \tilde u''$. By \eqref{eq:est-u''} and the fact that $[\tilde \lin,\psi']$ is supported on $[\ubar a, \ubar a+1]$, we have for $\tau$ small enough,
 	\begin{align}
	\notag\|  E_n\|_0 &\leq C e^{\gamma (\ubar a+1)} \|[\tilde \lin, \psi']  u'': C^{0,\alpha} (\tilde \Sigma, g_{\tilde \Sigma}) \| \\
	\notag	&\leq C e^{\gamma (\ubar a+1)}  \tau^{-1-\alpha}e^{- \gamma  (5\delta_s/\tau -1)} \| E_{n-1} \|_0 \\
	\label{eq:contraction}	& \leq e^{-\delta_s/\tau} \| E_{n-1}\|_0.
	\end{align}

We define $ v_E := \sum_{n=1}^{\infty}  v_n$ and $b_{E} := \sum_{n=1}^{\infty} b_n$. The three series converge and we have the desired estimates from \eqref{eq:contraction} and Proposition \ref{prop:7-1}. The function $ v_E$ is uniquely determined from the construction and satisfies $\tilde \lin  v_E  = E + \Theta(b_E)$. 
\end{proof}

\section{Quadratic Term}
\label{sec:quadratic}

\begin{proposition}
\label{prop:8-1}
Given $ v \in C_{loc}^{2,\alpha}(\tilde M)$ with $\|  v \|_2$ smaller than a suitable constant, the graph $\tilde M_v$ of $ v$ over $\tilde M$ is a smooth immersion; moreover
	\[
	\| \tilde H_v + \tilde X_v \cdot \nu_v - (\tilde H+ \tilde X \cdot \nu) - \tilde \lin   v\|_0 \leq C \| v\|_2^2,
	\]
where $\tilde H$ and $\tilde H_v$ are the mean curvature of $\tilde M$ and $\tilde M_v$ pulled back to $\tilde M$, respectively, and similarly, $\nu$ and $\nu_v$ are the oriented unit normal of $\tilde M$ and $\tilde M_v$ pulled back to $\tilde M$.
\end{proposition}

\begin{proof}
On the bounded piece $\tilde \Sigma \cup \tilde \ccal \cup \tilde \dcal$, the result follows from formulas for normal variations of $\tilde H$ and $\nu$ (see Appendix B \cite{kapouleas;embedded-minimal-surfaces} or  Section 4.2 \cite{mine;self-trans}). On the outer plane $\tilde \pcal$, a simple computation using $\xi$ as a coordinate on $\tilde \pcal$ shows that 
	\begin{align*}
	\tilde H_v + \tilde X_v \cdot \nu_v -  \tilde \lin   v & = \left(\delta_{ij} - \frac{D_{\xi_i}  v D_{\xi_j}  v}{1 + |D v|^2} \right) D^2_{\xi_i \xi_j} v - \xi \cdot D v +  v - (\Delta  v - \xi \cdot D v + v ) \\
	 	&= - \frac{D_{\xi_i} v D_{\xi_j}  v}{1 + |Dv|^2} D^2_{\xi_i \xi_j}  v.
	\end{align*}
Therefore, we have
	\[
	\| \tilde H_v + \tilde X_v \cdot \nu_v - (\tilde H+ \tilde X \cdot \nu) - \tilde \lin   v\|_0  \leq C \|  v \|_2^2.  \qedhere
	\]

\end{proof}

\section{The fixed point argument}
\label{sec:fpt}

We are now ready to prove the main result of this paper.

\begin{vartheorem} 

There exist a natural number $\bar  m$ and a constant $\zeta >0$  so that for any natural number $m>\bar m$, there exist   a constant  $b\in [-\zeta \sqrt 2/m, \zeta \sqrt 2/m]$ and  a smooth function $v$ on the initial surface $\tilde M(b, \sqrt 2/m)$ defined in Section \ref{ssec:initial-surface} such that the graph $\tilde M_m$ of $v$ over $\tilde M(b, \sqrt 2/m)$ has the following properties:
	\begin{enumerate}
	\item $\tilde M_m$ is a complete smooth surface which satisfies the equation $\tilde H + \tilde X \cdot \nu =0$.
	\item $\tilde M_m$ is invariant under rotation of $180^{\circ}$ around the $\tilde x$-axis.
	\item $\tilde M_m$ is invariant under reflections across planes containing the $\tilde z$-axis and forming angles $\pi/(2m) + k \pi/m$, $k \in \mathbf Z$, with the $\tilde x$-axis.
	\item Let $U =B_2 \cap \{\tilde z>0\}$ be the open top half of the ball of radius $2$. As $m \to \infty$, the sequence of surfaces $\tilde M_m$ tends to the  sphere of radius $\sqrt 2$ centered at the origin on any compact set of  $U$. 
	\item $\tilde M_m$ is asymptotic to a cone.
	\item If we denote by $T$ the translation by the vector $-\sqrt 2 \vec e_x$, the sequence of surfaces $m T(\tilde M_m) = \{ (m\tilde x, m \tilde y, m \tilde z) \mid (\tilde x + \sqrt 2 \vec e_x, \tilde y, \tilde z) \in \tilde M_m \} $  converges in $C^k$ to the original Scherk surface $\Sigma_0$ on compact sets, for all $k \in \mathbf N$.
	\end{enumerate}	
\end{vartheorem}

\begin{proof}
Let us denote by $\tau$ the quantity $\sqrt 2 /m$. We fix $\alpha' \in (0 ,\alpha)$ and define the Banach space
	\[
	\chi=C^{2, \alpha'} (\tilde M(0, \tau)).
	\]
Denote by $D_{b, \tau}: \tilde M(0, \tau) \to \tilde M(b, \tau)$ a family of smooth  diffeomorphisms which depend smoothly on $b$ and satisfy the following conditions: for every $f \in \mathcal C^{2, \alpha}(\tilde M(0, \tau))$ and $f' \in C^{2,\alpha}(\tilde M(b, \tau))$, we have
	\begin{align*}
	\| f \circ D_{b,  \tau} ^{-1} \|_2 \leq C\| f \|_2, \quad \| f' \circ D_{b,  \tau}\|_2 \leq C \|f'\|_2.
	\end{align*}
The diffeomorphisms $D_{b, \tau}$ are used to pull back  functions and norms from $\tilde M(b, \tau)$ to  $\tilde M(0, \tau)$. 

We fix $\tau$ and omit the dependence in $\tau$ in our notations of maps and surfaces from now on. Let
	\[
	\Xi = \{(b, u) \in \R\times \chi: |b| \leq \zeta  \tau, \| u \|_2 \leq \zeta \tau\},
	\]
where $\zeta$ is a large constant to be determined below. 
The map $\ical: \Xi \to \R  \times \chi$ is defined as follows. Given $(b,  u) \in \Xi$, let $v=u \circ D_{b}^{-1}$, $\tilde M = \tilde M(b)$ and let  $\tilde M_v$ be the graph of $v$ over $\tilde M$. We define the function $\F:\R  \times C^{2,\alpha}(\tilde M, g_{\tilde M}, e^{-\gamma s}) \to \R$ by 
	\[
	\F(b, v) = \tilde H_v + \tilde X_v \cdot \nu_v,
	\]
where $\tilde H_v$ and $\nu_v$ are the mean curvature and the oriented unit normal of $\tilde M_v$, respectively pulled back to $\tilde M$. Proposition \ref{prop:8-1} asserts that
	\[
	\|\F(b, v) - \F (b, 0) -  \tilde \lin_{\tilde M} v \|_0 \leq C \zeta^2  \tau^2.
	\]
Applying Theorem \ref{thm:linear-invertible} with $E = \F(b, v) - \F (b, 0) -  \tilde \lin_{\tilde M} v $, we obtain $v_E$ and $b_E$ such that 
	\begin{gather*}
	\tilde \lin_{\tilde M} v_E = E  + \Theta(b_E),\\
	\|v_E \|_2 \leq C \zeta^2  \tau^2, \quad |b_E | \leq C \zeta^2  \tau^2.
	\end{gather*}
Hence,
	\[
	\F(b,  v) = \F(b,  0) + \tilde \lin_{\tilde M} v + \tilde \lin_{\tilde M} v_E  - \Theta(b_E).
	\]
Propositions \ref{prop:smooth-dependence} and \ref{prop:4-20}, and Theorem \ref{thm:linear-invertible}  give us $v_H$ and $b_H$ satisfying 
	$
	\tilde \lin_{\tilde M} v_H = \F(b, 0) + \Theta(b_H), 
	$
 so
 	\[
	\F(b,  v) =  \tilde  \lin_{\tilde M} v + \tilde \lin_{\tilde M} v_H+ \tilde \lin_{\tilde M} v_E -\Theta(b_E + b_H).
	\]
We define the map $\mathcal I: \Xi \to \R \times \chi$ by 
	\[
	\mathcal I(b,  u) = (b-b_E - b_H, (-v_E-v_H)\circ D_{b} ).
	\] 
We now arrange for $\mathcal I (\Xi )  \subset \Xi$. Since
	\begin{align*}
	\| - v_E -v_H \|_2 &\leq C( \tau +  \zeta ^2  \tau^2),\\
	|b-b_E - b_H| &\leq  C( \tau +  \zeta ^2  \tau^2),
	\end{align*}
we can choose  $\zeta > 2C$ and $ \tau < \zeta^{-2}$ in order to get $C( \tau +  \zeta ^2  \tau^2) < \zeta  \tau$.

The set $\Xi$ is clearly convex. It is a compact set of $\R\times \mathcal X$ from the choice of the H\"older exponent $\alpha'<\alpha$. The map $\mathcal I$ is continuous by construction therefore we can apply the Schauder Fixed Point Theorem (p. 279 in \cite{gilbarg-trudinger;elliptic-equations}) to obtain a fixed point $(b_{ \tau},  u_{ \tau})$ of $\mathcal I$ for every $ \tau \in (0, \delta_{\tau})$ with $\delta_{\tau}$ small enough. The graph of $v=u_{ \tau} \circ D_{b,  \tau}^{-1}$ over the surface $\tilde M(b_{\tau},\tau)$ is  then a self-shrinking surface. It is a smooth surface by the regularity theory for elliptic equations. The properties (ii) and (iii) follow from the construction. 
\end{proof}


\bibliographystyle{siam}
\def\cprime{$'$}

\end{document}